\documentclass[11pt]{amsart}

\usepackage[T1]{fontenc}
\usepackage{amssymb,amsmath,amsthm,amsfonts,graphicx,hyperref,url,dsfont}

\usepackage{listings}
\usepackage{xcolor}
\usepackage{hyperref}
\usepackage[figure]{hypcap}

\usepackage{siunitx}
\usepackage{booktabs}
\usepackage[english]{babel}
\usepackage{blindtext}
\RequirePackage{fix-cm}

\usepackage{fullpage}

\usepackage{enumitem}

\newcommand{\eps}{\varepsilon}
\newcommand{\vphi}{\varphi}
\newcommand{\net}[1]{(#1)_{\varepsilon \in (0,1]}}

\theoremstyle{theorem}
\newtheorem{thm}{Theorem}[section]
\newtheorem{lem}[thm]{Lemma}

\newtheorem{de}{Definition}[section]

\newtheorem{cor}{Corollary}[section]

\newtheorem{ass}{Assumption}[section]

\theoremstyle{definition}
\newtheorem{rmk}{Remark}
\newtheorem{ex}{Example}[section]

\def\R{\mathbb R}

\def \E{\mathcal E}
\def\N{\mathbb N}
\def \F{\mathcal F}
\def \P{\mathbb P}
\def \B{\mathcal B}

\def\D{\mathcal D}

\def\tilde{\widetilde}

\def\endproof{\hfill$\Box$\\}
\def\supp{\mathop{\rm supp}\nolimits}    
    
\usepackage[T1]{fontenc}
\usepackage[utf8]{inputenc}
\usepackage{combelow}

\usepackage{newunicodechar}

\usepackage{blindtext}
\usepackage{float}

\newcommand{\normtri}[1]{{\left\vert\kern-0.25ex\left\vert\kern-0.25ex\left\vert #1
		\right\vert\kern-0.25ex\right\vert\kern-0.25ex\right\vert}}		

\newcommand{\norm}[2]{{\left\| #1 \right\|}_{#2}}

\newcommand{\ltwo}[1]{\norm{#1}{L^2}}

\newcommand{\linf}[1]{\norm{#1}{L^\infty}}

\begin{document}

\title[Stochastic very weak solution to parabolic equations with singular coefficients]{Stochastic very weak solution to parabolic equations with singular coefficients}

\author[S. Gordi\'c, T. Levajkovi\'c, Lj. Oparnica]{Sne\v zana Gordi\'c$^1$, Tijana Levajkovi\'c$^2$, \\ Ljubica Oparnica$^3$}
\date{\today}

\thanks{$^1$Sne\v zana GORDI\' C, {Faculty of Education, University of Novi Sad, Podgori\v cka 4, 25000 Sombor, Serbia, \texttt{snezana.gordic@uns.ac.rs, gordicsnezana@gmail.com}}\\
\indent$^2$Tijana LEVAJKOVI\'C, {Institute of Stochastics and Business Mathematics, Vienna University of Technology,  Austria, \texttt{tijana.levajkovic@tuwien.ac.at}}\\
\indent$^3$Ljubica OPARNICA,  \texttt{ljubica.oparnica@pef.uns.ac.rs,ljubica.oparnica@gmail.com} \& Faculty of Education, University of Novi Sad, Podgori\v cka 4, 25000 Sombor, Serbia
}

\maketitle


\begin{abstract}
A class of stochastic parabolic equations with singular potentials is analyzed within the chaos expansion framework, utilizing the Wick product to handle the multiplication of generalized stochastic processes. The analysis combines the chaos expansion method from white noise analysis with the concept of very weak solutions from partial differential equation theory. The stochastic very weak solution to the parabolic evolution problem is defined, and its existence and uniqueness are established. For sufficiently regular potentials and data, we demonstrate the consistency of the stochastic very weak solution with a stochastic weak solution. An illustrative example is provided, potential applications are reviewed, and future challenges are outlined.
\end{abstract}


\section{Introduction and preliminaries} \label{sec1} 

Stochastic processes are a fundamental concept in probability theory with broad applications across fields such as finance and physics. However, certain stochastic processes that emerge in applications cannot be defined in the classical sense. One notable example is white noise, which effectively models random fluctuations in dynamic systems.
White noise, interpreted as the time derivative of Brownian motion, cannot be defined as a classical process and instead, it is treated as a generalized stochastic process (GSP).
Drawing on the foundational theory of GSPs, particularly those developed via chaos expansions by Hida, Kondratiev, \O ksendal, Kuo, and others, this work focuses on a specific class of GSPs.
Specifically, we employ GSPs of Kondratiev-type to analyze the stochastic parabolic evolution problem:
\begin{equation*}
\left(\frac{\partial}{\partial t}  - {\mathcal L} \right) \, U  + \,Q\, \cdot \, U = F,\quad \quad 
U |_{t=0} = G,
\end{equation*}
where $U$ is an unknown GSP,  while the potential $Q$,  the driving force $F$ and the initial condition $G$ are given GSPs. The operator $\mathcal{L}$ represents an elliptic operator, and the product $\cdot$ is  interpreted based on the specific context of the problem.

Generalized stochastic processes singular with respect to both spatial and random arguments 
were first introduced in \cite{PS2007}, where they were defined as linear continuous mappings from a deterministic test function space into a Kondratiev space of generalized random variables. This framework has been applied in various contexts. For instance, in \cite{Sel2008}, the stochastic Dirichlet problem was studied, and an appropriate algebra of generalized stochastic processes was constructed. The fundamental solutions of the stochastic Helmholtz equation were provided in \cite{Sel2011}. 

Building on these foundational works, we focus on the analysis of a stochastic evolution initial value problem where the generalized stochastic potential  $Q$ is singular in space. Specifically, we study the stochastic evolution problem:
\begin{equation}\label{Eq:ISH3}
\begin{split}
\left(\frac{\partial}{\partial t}  - {\mathcal L} \right) \, U(t, x, \omega) + Q( x, \omega) \, 	\lozenge \, U(t, x, \omega) &= F(t, x, \omega), \\
U(0, x, \omega) &= G(x, \omega),
\end{split}
\end{equation}
where $t>0$ and $x\in \R^d$ are deterministic variables, $\omega\in\Omega$ is an elementary event  from the probability space $(\Omega, {\mathcal F}, P)$. The data $F$ and $G$ are  given GSPs  depending on space and/or time, $Q$  is given GSP singular in space, and $U$ is the unknown GSP. The operator $\mathcal{L}$ denotes an elliptic operator acting on the space domain, and $\lozenge$ denotes the Wick product. 

The motivation for studying this problem stems from earlier works. For 
 $Q$ bounded in space, the problem has been analyzed in \cite{GLO2023} using the Wick product to address the ill-defined multiplication of 
 $Q$ and $U$. For deterministic singular potentials, such as Dirac delta distributions, the problem was studied in \cite{GLO2021}, where the concept of a stochastic very weak solution was introduced to handle the ill-defined multiplication of distributions.

In this paper, we extend these analyses to cases where  $Q$ is both stochastic and singular in space. To address the associated challenges, we employ the method of chaos expansions from white noise analysis \cite{HKPS2013, HOUZ1996} and the framework of very weak solutions \cite{ARST2020c, GR15}. These tools enable a rigorous approach to the problem and provide a foundation for understanding the interaction of singular stochastic potentials with generalized stochastic processes.

Previous works on similar problems that using different methodologies include the use of pseudodifferential and Fourier integral operators, tools derived from microlocal analysis, in \cite{Coriasco2021} to establish the existence of solutions for a class of linear parabolic SPDEs with polynomially bounded variable coefficients. Similarly, in \cite{AWZ1998}, the existence of unique solutions to parabolic SPDEs driven by Poisson white noise was proven by interpreting these equations as stochastic integral equations of the jump type involving evolution kernels. The asymptotic behavior as  $t \to \infty$  of solutions to heat equations with white noise potentials was investigated in \cite{hu2002}.

Equation \eqref{Eq:ISH3}  arises in the modeling of various phenomena, including disordered media with random impurities, environmental transport with random sources, biological tissues with random heterogeneities, quantum physics (e.g., Anderson localization), reaction-diffusion systems with random catalysis, and random risk potentials. These applications span a wide range of fields such as biophysics, neuronal dynamics, quantum physics, groundwater contamination, and economic and financial modeling (see Section 6).

The  paper is organized as follows. We close this section by fixing  notations, recalling the necessary definitions, and sumerising previous results. In Section \ref{sec:boundedQ}, we review and improve the results on the problem \eqref{Eq:ISH3} with random and bounded in space potential $Q$, given in  \cite{GLO2023}.  In Section \ref{sec:singularQ_DefExi}, we introduce the regularisation nets of generalized stochastic processes, define stochastic very weak solution for the problem \eqref{Eq:ISH3}, and show its existence. Section \ref{sec:singularQ_Uniqueness} contains the uniqueness results of the stochastic very weak solution in an appropriate sense. In Section \ref{sec:singularQ_Consistency}, we show the consistency of the stochastic very weak solution and the stochastic weak solution for problem \eqref{Eq:ISH3} if the  potential $Q$ is bounded in space. 
In Section \ref{sec6}, we provide an example by analyzing equation \eqref{Eq:ISH3} with the Laplacian as the operator $\mathcal{L}$, 
and we discuss various potential applications across diverse fields demonstrating its wide-ranging relevance. Finally, we  highlight potential directions for future research. 

\subsection{Notations and classical results}
 For the later use, and completeness of the presentation of results, we fix well-known notation, and recall important notions. As usual, by $L^2(\R^d)$ we denote the space of square integrable functions over $\R^d$, by  $L^\infty(\R^d)$ we denote the space of essentially bounded measurable functions, 
 by $\mathcal{D}(\R^d):=\mathcal{C}_{0}^{\infty} (\mathbb{R}^d)$ we denote the space of compactly supported smooth functions, the space of distributions we denote by $\mathcal{D}'(\mathbb{R}^d)$, by $\mathcal{S} (\mathbb{R}^d)$ and $\mathcal{S}' (\mathbb{R}^d)$ we denote  the Schwartz space of rapidly decreasing functions and the Schwarz space of tempered distributions, respectively, and  by $\mathcal{E}'(\mathbb{R}^d)$ we denote the space of compactly supported distributions. 
 
 Further, we will use the vector valued functions and distributions. 
The space  $L^2([0,T]; L^2(\mathbb R^d))$ denotes  the Hilbert space of square integrable functions over $[0,T]$  with values in $L^2(\mathbb R^d)$, and with the norm 
$$\|f\|_{L^2([0,T]; L^2(\mathbb R^d))}=   \left( \int_0^T \|f(t, \cdot)\|^2_{L^2(\mathbb R^d)} \, dt \right)^{1/2}
=  \left( \int_{\mathbb R^d } \, \int_0^T   |f(t, x) |^2  \, dt \, dx \right)^{1/2}.$$ 
The space $C([0,T];  L^2(\mathbb R^d))$ denotes the Banach space of continuous functions over $[0,T]$  with values in $L^2(\mathbb R^d)$ with the norm  $\|f\|_{C([0,T]; L^2(\mathbb R^d))} = \sup_{t \in [0,T]} \|f(t, \cdot)\|_{L^2(\mathbb R^d)},$ $L^1(0,T; {L^{2}(\mathbb R^d)})$ is the space of integrable functions over $[0,T]$ with values in  $L^{2}(\mathbb R^d),$ and $AC([0,T],L^2(\R^d))$ is the space of absolutely continuous functions over $[0,T]$, i.e., 
differentiable functions in $t$, a.e. on $[0,T]$,  with integrable first derivative,  and with values in $L^2(\R^d)$. 

For $m \in \N$ and $1 \leq p \leq +\infty$, we denote by $H^{p,m}(\R^d)$ the set of functions $u$ on $\R^d$ such that all derivatives of $u$ up to the order $m$  belong to $L^p (\R^d),$ and the space $H^{p,m}(\R^d)$ is equipped with
	 $$\|u\|_{p,m}=\left( \sum_{|\alpha| \leq m} |(\partial/\partial x)^\alpha u(x) |^p \, dx\right)^{1/p}.$$
Note that $\mathcal{D}(\R^d)$ is dense in $H^{p,m}(\R^d).$ We denote by $H^{m} (\R^d)$ the Sobolev space $H^{2,m} (\R^d).$
For $s>0$ real, the Sobolev space $H^{-s} (\R^d)$ is the space of tempered distributions $u$ whose Fourier transform $\hat u$ is a square integrable function with respect to measure $(1+|\xi|^2)^{s/2}\,d\xi$. $H^{-s} (\R^d)$ is a Hilbert space for all real $s>0$. Let $H^{-s}_c (\R^d)$ be the spaces of distributions $u \in  H^{-s} (\R^d)$ such that $\supp u$ is compact in $\R^d.$ If $u \in H^{-s} (\R^d)$ there are functions  $k_{\alpha} \in L^2(\R^d)$ such that 
\begin{equation} \label{rep_u}
u= \sum_{|\alpha|\leq s} ({\partial}/{\partial x})^\alpha k_{\alpha}.
\end{equation} 
On $H^{-s} (\R^d)$ the following norm is defined
 \begin{equation} \label{H} 
 \| u\|_{H^{-s}_c}^2 = \inf \sum_{|\alpha|\leq  s} \| k_{\alpha}\|_{L^2}^2,
 \end{equation}
 where the infimum is computed over all the representations of $u$  of the form \eqref{rep_u}.
$H^{-s}(\R^d)$ is the dual space of the Sobolev space $H^s (\R^d)$. Note that $\E'(\R^d) = \bigcup_{s \in \R} H^{-s}_c (\R^d).$  

%


To establish a solid foundation for the analysis, we will outline key assumptions that underpin our model under the special environment Assumption. Each of the assumptions will be clearly stated ensuring coherence and relevance to the context of the study. We start with the assumption on the  operator $ \mathcal L $.

\begin{ass} \label{op_L}
	The operator $ \mathcal L $ is unbounded and closed operator on  $L^2(\R^d)$ with dense domain $D\subseteq L^2(\mathbb R^d)$,  which generates  a $C_0$-semigroup $(T_t)_{t\geq 0}$ on $L^2(\R^d)$, and  $w \in \mathbb{R}$ and $M>0$ are  the stability constants from the semigroup estimate
	\begin{equation} \label{SemiGrpEst}
		\|T_t\|_{\mathfrak{L}(L^2(\R^d))}\leq Me^{wt}, \quad t\geq 0.
	\end{equation}
\end{ass}
%
 The well-known example of the operator $\mathcal{L}$ is the Laplace operator $\Delta$ on $L^2(\R^d)$, with dense domain $D=H_0^1 (\R^d)$ which generates a heat $C_0$-semigroup of contractions with the stability constants $M=1$ and $w=0$. In this case, equation \eqref{Eq:ISH3} becomes the stochastic heat equation.  

We recall the classical result from \cite{Pazy1983}, which establishes the existence of a unique solution for parabolic initial value problem  obtained via semigroup approach, and we use notions and estimates from \cite[Theorem 3]{GLO2021}.

\begin{thm}\label{T:DetEq}
Given the force term $f\in AC([0,T];L^2(\R^d))$, 
the initial condition $g\in D$  and the potential $q \in L^{\infty} (\mathbb{R}^d)$,  
the (deterministic) parabolic initial value problem
	\begin{equation*} 
	\left( \frac{\partial}{\partial t}  -  {\mathcal L} \right)  u(t, x) + q (x)  u(t, x) = f(t, x ), 
	\quad u(0, x) = g(x),
	\end{equation*}
	has a unique bounded solution $u\in AC([0,T];L^2(\R^d))$ 
	which satisfies  
	\begin{equation*} 
	\|u(t,\cdot)\|_{L^2} \leq  M(t)  \left( \|g(\cdot)\|_{L^{2}} +\int_0^t \|f(s,\cdot)\|_{L^2} \,ds\right), \quad t\in (0,T],
	\end{equation*}
	where 
\begin{equation}\label{eModt}
	M(t) : = M \exp {\left( \left(w + M \|q\|_{L^\infty}\right)t \right)}, \quad t\in (0,T].
\end{equation}
\end{thm}
Further, if we denote
\begin{equation}\label{mtilda}
\tilde{M}(t) := \displaystyle \int_0^t M(s) \,ds = \tilde{M}(t) = \dfrac{M(t)-M}{w+M\|q \|_{L^\infty}},  \qquad \tilde{M} = \tilde{M} (M,w,\|q\|_{L^\infty}),
\end{equation}
then the following  holds.
\begin{lem} \label{Lema M(t)} 
If $M(t)$ is as in \eqref{eModt}, 
 and  $\tilde{M}(t)$ as in \eqref{mtilda}, then the following estimates hold:
 
		$\displaystyle \int_0^t  M(s)  \tilde{M}(s)^n  \,ds \leq	 \  \tilde{M}(t)^{n+1}$ and   $\displaystyle \int_0^t s M(s) \tilde{M}(s)^n\,ds \leq  t \tilde{M}(t)^{n+1}$ for all $n \in \N_0.$
\end{lem}

\subsection{White noise analysis}
In this section, we briefly revisit fundamental concepts from white noise analysis, and for  details we refer  to \cite{HOUZ1996}.
\subsubsection{Multi-indices of arbitrary lenght}
Define $\mathcal{I}$ as the set of sequences of non-negative integers which have finitely many nonzero components, i.e., 
$\gamma= (\gamma_1,\gamma_2, \ldots,\gamma_m,0,0,\ldots),$ $\gamma_i \in \N_0$, $i=1, 2,\ldots, m,$ and $\gamma_i=0$ for $i>m,$ $m\in\N.$ Thus, $\mathcal{I}:=(\N_0^{\N})_c $, $\N_0=\N \cup \{0\}$, with $(\cdot)_c$ denoting compactly supported sequences.
The zero multi-index, denoted by $\textbf{0}$, has all coordinates equal to 0, and the $k$th unit multi-index is a multi-index whose coordinates are all zero, except  the $k$th coordinate, which is equal to 1, and it is denoted by $e_k$,  $k\in \N$.  The length of multi-index $\gamma \in \mathcal{I}$  is  $|\gamma|= \sum_{i=1}^{\infty}\gamma_i$,  and the factorial of $\gamma$ is $\gamma ! = \prod_{i=1}^{\infty} \gamma_i !$.
As in the case of usual multi-indices the  Cauchy–Schwarz inequality  holds, i.e.,
	\begin{equation}\label{CSmulti}
		\left(\sum_{\alpha \in \mathcal{I}} |x_\alpha y_\alpha| \right)^2 \leq \left(\sum_{\alpha \in \mathcal{I}} |x_\alpha|^2\right) \left(\sum_{\alpha \in \mathcal{I}} |y_\alpha|^2\right), \quad x_\alpha, y_\alpha \in \R.
	\end{equation}
 For $\alpha=(\alpha_1, \alpha_2,\dots)\in \mathcal{I}$ and $\beta=(\beta_1,\beta_2,\dots)\in \mathcal{I}$, \textcolor{black}{we write that $\alpha \leq \beta$ if}  $\alpha_i \leq \beta_i$ for all $i \in \N$. For $x=(x_1,x_2,\dots,x_m) \in \R^m$ and $m$-dimensional multi-index $\gamma = (\gamma_1,\gamma_2, \dots,\gamma_m)$, we have $x^\gamma = x_1^{\gamma_1} \cdot x_2^{\gamma_2} \cdot \dots \cdot x_m^{\gamma_m}.$ If $c$ is a real positive constant and $\gamma
	\in \mathcal{I}$, then $c^\gamma = \prod_{k=1}^\infty c^{\gamma_k}.$

For $\gamma\in\mathcal{I}$ one defines the weight $(2\N)^\gamma$ as
$$(2\N)^\gamma := \prod_{k=1}^\infty (2k)^{\gamma_k}.$$
Note that, the above product and the sums in the  Cauchy–Schwarz inequality are finite, since all  multi-indices $\gamma \in \mathcal{I}$  are of the finite lenght. The following lemma collects results on multi-indices of arbitrary length, and the weight $(2\N)^\gamma$. 
\begin{lem}\label{L:mi prop} 
\textcolor{black}{The following hold:}
	\begin{itemize}
	\item[a)]  \textcolor{black}{For all $\gamma \in \mathcal{I},$} $|\gamma|\leq (2\N)^{\gamma}.$ 
	\item[b)] For every $c>0$, there exists $s \geq 0$ such that $c^{|\gamma|} \leq (2\N)^{s\gamma}$ and $c^\gamma \leq (2\N)^{s\gamma}$ for \textcolor{violet}{all} $\gamma \in \mathcal{I}.$
	\end{itemize}
\end{lem}
\begin{rmk}\label{Rem_on_s}
Let us give a note on the choice of $s$ in Lemma \ref{L:mi prop} b). 
 For arbitrary $c>0$, one has to find $s = s(c)$ such that $c^{|\gamma|}\leq (2\N)^{s\gamma}$.
To find such $s$ note that there is a $k\in\N$ such that $2k<c<2(k+1)$. For $i \leq k$ we look for 
$s_i$ so that $c^{\gamma_i}\leq (2i)^{s_i\gamma_i}$ which is the same as $s_i\geq\frac{\ln c}{\ln (2i)}$. Setting  $s_i = 1$ for $i>k$ and choosing 
\begin{equation}\label{s}
s = \frac{\ln c}{\ln 2}+1
\end{equation}
we have that $s > s_i$ for all $i\in\N$ and therefore 
$$ c^{|\gamma|} = c^{\gamma_1} c^{\gamma_2} \cdots c^{\gamma_n} \leq 
\prod_{i=1}^{n} (2i)^{s_i\gamma_i} \leq (2\N)^{s\gamma}. $$
\end{rmk}
We recall the following result from \cite[Proposition 2.3.3, p. 35]{HOUZ1996}.
\begin{lem} \label{L:konvCp}
		$\sum_{\gamma \in \mathcal{I}} (2\N)^{-p\gamma}< \infty$ if and only if $p>1.$ 	
\end{lem}

In our work, the constant from the previous lemma
\begin{equation}\label{Cp}
	C_p:= \sum_{\gamma \in \mathcal I} (2\N)^{-p\gamma}, \quad p>1,
\end{equation}
plays an important role. Since $C_p \to 0$ as $p \to \infty$, there exist $\delta>0$ and $p_0 \in \N$ such that $C_p <\delta$ for all $p>p_0.$ In particular, the following lemma will be commonly used.
\begin{lem} \label{lemanew}
Let $d(\gamma)=  |\gamma|^m d^{n\gamma}$ or $d(\gamma)=  |\gamma|^m d^{n|\gamma|}$, where $ n,m \in \mathbb{N}$, and $d>0$ is a constant. Define $$D_p:= \sum_{\gamma \in \mathcal I}  d(\gamma) (2\N)^{-p\gamma}.$$ Then, there exists $s\geq 0$ such that $$D_p \leq C_{p-m-sn}.$$
\end{lem}

\proof Let $d(\gamma)=  |\gamma|^m d^{n\gamma}$. By Lemma \ref{L:mi prop} b) we obtain that there exists $s\geq 0$ such that
\begin{equation*}
D_p=	\sum_{\gamma \in \mathcal I} |\gamma|^m d^{n\gamma} (2\N)^{-p\gamma}\leq 
	\sum_{\gamma \in \mathcal I}(2\N)^{m\gamma} (2\N)^{sn\gamma} (2\N)^{-p\gamma} = 
	\sum_{\gamma \in \mathcal I}(2\N)^{-(p-m-sn)\gamma} = C_{p-m-sn},
\end{equation*}
which is finite for $p-m-sn>1$, i.e., $p>1+m+sn$, and $s=\frac{\ln d^n}{\ln 2}+1$ by \eqref{s}. The proof is similar for  $d(\gamma)=  |\gamma|^m d^{n|\gamma|}$.
\endproof

We note that, for example, when $m=2$, from Lemma \ref{lemanew}, we have $$D_p= \sum_{\gamma \in \mathcal I} |\gamma|^2 d^{n\gamma} (2\N)^{-p\gamma}\leq C_{p-2-sn},$$ and this sum is finite for $p>3+sn,$ with $s=\frac{\ln d^n}{\ln 2}+1$. Similarly, for the case when $m=0$ and $n=0$, we obtain $D_p=C_p$, which is finite for $p>1,$ as it follows from Lemma \ref{L:konvCp}.

\subsubsection{Chaos expansions in $L^2(\mu)$}

%

The Gaussian white noise probability space,  $(\Omega,\F,\P)$, is  identified  with  $(\mathcal{S}'(\R^d),\B,\mu)$, where  $\B$ is the family of all Borel subsets of $\mathcal{S}'(\R^d)$ equipped with the weak$^*$ topology and $\mu$ is the Gaussian measure on $\mathcal{S}'(\R^d)$ corresponding to the characteristic function $$\int_{\mathcal{S}'(\R^d)} e^{i\langle \omega, \phi \rangle} \,d\mu(\omega)= \exp \left[ -\frac{1}{2} \| \phi \|^2_{L^2(\R^d)}\right], \quad \phi \in \mathcal{S}(\R^d),$$ whose existence is guarantied by the Bochner--Minlos theorem \cite{HOUZ1996}. 

The  space of square integrable random variables on  $(\Omega,\F,\P)$, denoted by $L^2(\mu)$,   is a Hilbert space with the norm $\|\cdot\|_{L^2(\mu)}$ induced by the inner product 
$$(f, g)_{L^2(\mu)} = \mathbb{E} \, (fg) = \int_{\Omega} f(\omega) \, g(\omega) \, d\mu(\omega) <\infty, \quad \text{for}\,\,  f, g\in L^2(\mu),$$
where $\mathbb{E}$ denotes the expectation with respect to the Gaussian measure $\mu$. \\

For a given $\gamma\in \mathcal{I}$, the $\gamma$th {\it Fourier-Hermite polynomial} $H_\gamma$ is defined by
\begin{equation}\label{Halpha}
	H_{\gamma}(\omega)=\prod_{k=1}^\infty h_{\gamma_k}(\langle\omega,
	\eta_k\rangle), \quad\gamma\in\mathcal I, \quad \omega \in \Omega,
\end{equation}
where $\langle \cdot,\cdot \rangle$ denotes the action of $\omega \in \Omega= \mathcal{S}'(\R^d)$ on $\eta_k \in \mathcal{S}(\R^d)$, $\eta_k:=\xi_{\delta^{(k)}} = \xi_{\delta_1^{(k)}} \otimes \cdots \otimes \xi_{\delta_d^{(k)}}$ for $\delta^{(k)} = (\delta_1^{(k)}, \delta_2^{(k)}, \dots, \delta_d^{(k)})$ denotes the $k$th multi-index number in some fixed ordering of all $d$-dimensional multi-indices $\delta = (\delta_1, \dots, \delta_d) \in \N^d,$
$\xi_{\delta_j^{(k)}}$, $j=1,\dots,d$, are the Hermite functions of order $\delta_j^{(k)}$, and  $h_{\gamma_k}$, $k \in \mathbb N$, are the Hermite polynomials. 

The family $\{H_\gamma\}_{\gamma\in \mathcal{I}}$ forms an orthogonal basis in  $L^2(\mu)$, and $\|H_\gamma\|^2_{L^2(\mu)} = (H_\gamma,H_\gamma)_{L^2(\mu)}=\gamma!.$
\textit{The  Wiener-It\^o chaos expansion theorem} states that each square integrable random variable  $f\in L^2(\mu)$ has a unique representation of the form
\begin{equation*}
	f(\omega)  = \sum_{\gamma \in \mathcal{I}} \, c_{\gamma} \, H_{\gamma}(\omega),
	\enspace  \, \omega\in \Omega, \, c_\gamma = \dfrac{1}{\gamma!} (f,H_\gamma)_{L^2(\mu)}\in \mathbb{R}^d, \textcolor{violet}{ \text{}}
\end{equation*}
and that the isometry $\|f\|^2_{L^2(\mu)}=\sum_{ \gamma\in \mathcal{I}} \, c_{\gamma}^2 \, \gamma!$ holds. 
From the previous it follows  that $f$ with chaos expansion form  $	f(\omega)  = \sum_{\gamma \in \mathcal{I}} \, c_{\gamma} \, H_{\gamma}(\omega)$ belongs to $L^2(\mu)$ if and only if the growth condition $\sum_{ \gamma\in \mathcal{I}} \, c_{\gamma}^2 \, \gamma!<\infty$ is satisfied. Note that the expectation of the random variable $f$ is equal to the zeroth coefficient $c_\mathbf{0}$ in its chaos expansion form, i.e., $\mathbb{E}(f)=c_{\mathbf{0}}.$

\begin{ex} \label{ex1}
	Let $a := (a_k)_{k \in \N} \in l^2(\R)$, i.e., 
	$\sum_{k=1}^\infty a_k^2 < \infty$. Let $f=\sum_{k=1}^\infty a_k H_{e_k}$ be a Gaussian random variable, i.e., a random variable whose  chaos expansion representation has only multi-indices of length zero and one. It is an element in $L^2 (\mu)$ since $$\|f\|^2_{L^2(\mu)} = \sum_{k=1}^\infty a_k^2 \cdot (e_k)! = \sum_{k=1}^\infty a_k^2 < \infty.$$ 
\end{ex}

\subsubsection{Kondratiev spaces}
{The Kondratiev space of random variables} $(S)_{1,p}$, $p \in \N$, consists of the elements 
$f \in L^2(\mu)$ with the chaos expansion form $f = \sum_{\gamma \in \mathcal{I}} c_\gamma H_\gamma$, $c_\gamma \in \R^d,$ for which the following growth condition is satisfied
\begin{equation} \label{S1,p est}
	\Vert f \Vert_{p}^{2} := \sum_{\gamma \in \mathcal{I}}^{} c_{\gamma}^{2}(\gamma !)^{2} (2\N)^{\gamma p} < \infty.
\end{equation}
The space $(S)_{1,p}$ is Banach space with respect to the norm defined in \eqref{S1,p est}.
The projective limit of the spaces $(S)_{1, p}$ is the Kontratiev space $(S)_1$, i.e.,  
$(S)_1 \textcolor{violet}{:}= \varprojlim_{p\in\N} (S)_{1,p}$. 
A topology on $(S)_1$ is generated by the family of norms $\Vert \cdot \Vert _{ p}$, $p\in\N$, and $(S)_1$ is a nuclear space, \cite[Lemma 2.8.2]{HOUZ1996}. Thus, the Kondratiev space $(S)_1$  consists of the elements  $f \in L^2(\mu)$ with the chaos expansion form $f = \sum_{\gamma \in \mathcal{I}} c_\gamma H_\gamma$, $c_\gamma \in \R^d,$ such that \eqref{S1,p est} holds for all $p \in \N.$

The Kondratiev space of generalized random variables or stochastic distributions $(S)_{-1,-p},$ $p\in \N$, consists of random variables $F$ with
chaos expansions  $F=\sum_{\gamma \in \mathcal{I}} b_\gamma H_\gamma,$ $b_\gamma \in \R^d,$ such that
\begin{equation} \label{S-1-p est}
	\Vert F\Vert_{-p}^{2}: = \sum_{\gamma \in \mathcal{I}} b_{\gamma}^{2} (2\N)^{-\gamma p} < \infty. 
\end{equation}

The inductive  limit of the spaces $(S)_{-1,-p}$ is the Kontratiev space of distributions $(S)_{-1}$, i.e.,  
$(S)_{-1} = \varinjlim_{p \in \N}(S)_{-1,-p}$. It is the dual space of $(S)_1$, it is a nuclear Frech\'et space, and 
it consists of random variables with chaos expansion $F=\sum_{\gamma \in \mathcal{I}} b_\gamma H_\gamma,$ $b_\gamma \in \R^d,$ such that \eqref{S-1-p est} holds for some $p \in \N$.
The action of $F\in(S)_{-1}$ on $f\in(S)_1$ is given by
$
\langle F, f \rangle = \sum_{\gamma \in \mathcal{I}}  (b_\gamma,  c_\gamma) \gamma !,
$ where $(\cdot,\cdot)$ denotes the usual inner product in $\R^d.$
The spaces  $(S)_{1,p}$ and $(S)_{-1,-p}$ are separable  Hilbert spaces,  the embedding $(S)_{1,q} \subseteq (S)_{1,p}$ for $p\leq q$  is Hilbert-Schmidt type, and it holds 
$$
(S)_{ 1} \subseteq (S)_{1,q} \subseteq (S)_{1,p}  \subseteq L^2(\mu) \subseteq   (S)_{-1,-p}  \subseteq   (S)_{-1,-q} \subseteq   (S)_{ - 1}. 
$$ 

%
%
	
\subsubsection{$X$-valued generalized stochastic processes of Kondratiev-type}
Let $X$ be a normed vector space, and denote by $X\hat{\otimes}_\pi  L^2(\mu)$ the space of   square integrable stochastic $X$-valued processes, i.e., 
$$X \hat{\otimes}_\pi  L^2(\mu) =\left\{ F =\sum_{\gamma \in \mathcal{I}} f_\gamma H_\gamma \text{ with } f_\gamma \in X: \|f\|^2_{X \hat{\otimes}_\pi  L^2(\mu)} = \sum_{\gamma \in \mathcal{I}} \gamma! \|f_{\gamma}\|^2_{X}< \infty \right\},$$ 
where $\hat{\otimes}_\pi$ denotes the completion of the tensor product with respect to $\pi$-topology. For more details on tensor products and $\pi$-topology we refer to \cite{Traves}.

The space of $X$-valued stochastic processes of Kondratiev-type $X \hat{\otimes}_\pi (S)_{1,p}$, $p \in \N,$ is defined as
$$X \hat{\otimes}_\pi (S)_{1,p} \textcolor{violet}{:}= \left\{ F = \sum_{\gamma \in \mathcal{I}} f_\gamma H_\gamma \in X \hat{\otimes}_\pi L^2(\mu)
f_\gamma \in X
: \normtri{F}^2_{p} = \sum_{\gamma \in \mathcal{I}} \|f_\gamma\|^2_X \gamma!^2(2\N)^{p\gamma} < \infty \right\}, 
$$
and the space of $X$-valued test stochastic processes of Kondratiev-type $X \hat{\otimes}_\pi (S)_{1}$ {is defined} as the projective limit of the  tensor product spaces $X \hat{\otimes}_\pi (S)_{1,p }$, $p\in \N$, i.e., $$X \hat{\otimes}_\pi (S)_{1}\textcolor{violet}{:}= 
\varprojlim_{p\in\N}  X \hat{\otimes}_\pi (S)_{1,p}.$$

The space of $X$-valued generalized stochastic processes (GSP) of Kontratiev-type $X \hat{\otimes}_\pi (S)_{-1,-p},$ $p\in \N$, is defined as
 $$X \hat{\otimes}_\pi (S)_{-1,-p} \textcolor{violet}{:} = \left\{ F = \sum_{\gamma \in \mathcal{I}} f_\gamma H_\gamma: \normtri{F}^2_{-p} = \sum_{\gamma \in \mathcal{I}} \|f_\gamma\|^2_X (2\N)^{-p\gamma} < \infty \right\},
 $$ and 
 the space of $X$-valued GSP of Kondratiev-type 
   $X \hat{\otimes}_\pi (S)_{-1}$ {is defined} as the inductive limit of the  tensor product spaces $X\hat{\otimes}_\pi (S)_{-1,-p }$, $p\in \N$, i.e.,  
   $$X \hat{\otimes}_\pi (S)_{-1} \textcolor{violet}{:}= 
   \varinjlim_{p \in \N} X \hat{\otimes}_\pi (S)_{-1,-p}.$$

 %
 If  $F \in X \hat{\otimes}_\pi (S)_{-1}$, then $F \in X \hat{\otimes}_\pi (S)_{-1,-p_0}$ for some $p_0\in \N,$  and
\begin{equation}  \label{KondSN}
		\normtri{F}^{2}_{-p_0} = \sum_{ \gamma\in \mathcal{I}}\|f_\gamma\|^2_X (2\N)^{-p_0\gamma} < \infty,
\end{equation}
where $F = \sum_{ \gamma\in \mathcal I}f_\gamma H_\gamma$,  $f_\gamma\in X$.
If \eqref{KondSN} is finite for $p_0 \in \N$, then it is finite for all $p\geq p_0$. The minimal  $p_0$ such that \eqref{KondSN} holds is called {\em the critical exponent}. We denote the critical exponent of GSP of Kondratiev-type $F$ by $p_F$.

For a Banach space $X$, and for every $p \in \N$, the space $X \hat{\otimes}_\pi (S)_{-1,-p}$  is a Banach space, and the space $X \hat{\otimes}_\pi (S)_{-1}$ is a Frech\'et space. 
%
%
\textcolor{black}{Two $X$-valued (generalized) stochastic processes of Kondratiev-type $F_1$ and $F_2$ are equal if and only if 
the corresponding coefficients in their chaos expansions are equal.}


\begin{ex} \label{Ex:wn}
\begin{itemize}
\item[(a)] Let  $F \in AC([0,T]; L^2 (\R^d)) \hat{\otimes}_\pi (S)_{-1}$ be an $AC([0,T]; L^2 (\R^d))$-valued GSP of Kondratiev-type  with a chaos expansion form  $F(t,x,\omega) =\sum_{\gamma \in \mathcal I} f_\gamma (t,x)H_\gamma(\omega)$, $t>0$, $x \in \R^d$, $\omega \in \Omega = \mathcal{S}'(\R^d)$ with coefficients  $f_\gamma \in AC([0,T]; L^2 (\R^d))$, $\gamma\in \mathcal I$, and for some $p \in \N$ it holds
$$\normtri{F}^2_{-p}  = \sum_{\gamma \in \mathcal I} \|f_\gamma\|^2_{AC([0,T]; L^2 (\R^d))} (2\N)^{-p\gamma} < \infty.$$  Since $AC([0,T]; L^2 (\R^d))$ is a Banach space, it follows that $AC([0,T]; L^2 (\R^d)) \hat{\otimes}_\pi (S)_{-p}$ is a Banach space with  norm  $\normtri{\cdot}_{-p}$ and thus $AC([0,T]; L^2 (\R^d)) \hat{\otimes}_\pi (S)_{-1}$ is a Frech\'et space.
\item[(b)] The time  dependent white noise $W_t$ is an element of the space 
 $\mathcal{S}([0,T]) \hat{\otimes}_\pi (S)_{-1}$ with the chaos expansion $W_t(\omega) = \sum_{k=1}^\infty \xi_k (t) H_{e_k} (\omega),$ while the space dependent white noise is an element of the space $\mathcal{S}(\mathbb R^d) \hat{\otimes}_\pi (S)_{-1}$ with the chaos expansion $W_x (\omega) = \sum_{k =1}^\infty \eta_{k}(x)H_{e_k}(\omega)$,  
 where $\xi_k, \eta_k$ are as in the definition of the Fourier-Hermite polynomial $H_\gamma$ given by \eqref{Halpha}. In particular, for $d=1$, the space dependent white noise is given by $W_x(\omega) = \sum_{k=1}^\infty \xi_k (x) H_{e_k} (\omega)$.

\end{itemize}
\end{ex}


\subsubsection{Singular GSPs of Kondratiev-type} 

%
%
If $X$ is a Banach (or a Frech\'et) space, then by the nuclearity of a space $(S)_{1,p}$  (\cite[Lemma 2.8.2]{HOUZ1996}) it follows $$(X \hat{\otimes}_\pi (S)_{1,p})' \cong X' \hat{\otimes}_\pi (S)_{-1,-p},$$ where $X'$ denote the dual space of $X$. The space $X' \hat{\otimes}_\pi (S)_{-1,-p}$ is isomorphic to the space of linear bounded mappings $X \to (S)_{-1,-p}$, and it is isomorphic to the space $ \mathfrak{L}(X,(S )_{-1})$ of linear bounded mappings $(S)_{1,p} \to X'.$


\begin{de} \label{singGSP}  Let $X$ be a Frech\'et space.
	A singular GSP of Kondratiev-type 
is linear continuous mapping $X \to (S)_{-1}$, i.e., an element of the space $ \mathfrak{L}(X,(S )_{-1}).$
\end{de} 

Since $(S)_{-1}$ is a nuclear space, the space  $ \mathfrak{L}(X,(S )_{-1})$ is isomorphic to the space $X' \hat{\otimes}_\pi (S)_{-1}.$
Singular GSP of Kondratiev-type $\mathcal{E}'(\R^d) \hat{\otimes}_\pi (S )_{-1}$ and $\mathcal{S}'(\R^d) \hat{\otimes}_\pi (S)_{-1}$ are introduced in \cite{Sel2011}.
The process $F$ from $X' \hat{\otimes}_\pi (S)_{-1}$ has the chaos expansion from
$F = \sum_{\gamma \in \mathcal{I}} f_\gamma  H_\gamma$, with coefficients $f_\gamma \in X'$, $\gamma \in \mathcal{I}$,   satisfying	
\begin{equation}\label{seminorma_s_Bp}
\normtri{F}_{-p}^2 =\sum_{\gamma \in \mathcal{I}} \|f_{\gamma}\|^2_{X'} (2\N)^{-p\gamma} < \infty
\end{equation}
for some $p\in \N$, where $\| \cdot \|_{X'}$ denotes an appropriate seminorm in the space $X'.$ 

%
\begin{ex} \cite{Sel2011}
Let $K$ be a compact subset of $\R^d$ and let $\E_{m,K}(\R^d)$ be the space $C^{\infty}(\R^d)$ endowed with the topology given by the family of seminorms 
$$
p_{K,m} (f) =\displaystyle \sum_{|k| \leq m} \sup_{x \in K} |D^k f(x)|, \quad m \in \N, \enspace f \in C^{\infty}(\R^d).
$$
	%
	%
Let $\{K_k\}_{k \in \N}$ be the family of compact sets exhausting $\R^d$. The space $\E_m(\R^d)$ is the projective limit of  $\E_{m,K_k}(\R^d)$, $k \in \N$, 
 its dual space is  $ \E_m'(\R^d)$, i.e.,
$F \in \E_m'(\R^d)$  if there exist $K \subset \subset \R^d$ and $C>0$, such that $|\langle F, \varphi \rangle | \leq Cp_{K,m}(\varphi),$ for all $\varphi \in \E_m(\R^d).$
Let
\begin{equation*} 
\|F\|_{\E'_m} = \inf\{C>0: |\langle F, \varphi \rangle | \leq Cp_{K,m}(\varphi), \quad \text{for some compact set } K \subset \R^d\}.
\end{equation*}
The dual space $\E'(\R^d)$ of $\E(\R^d)$ is isomorphic to the inductive limit of the spaces $\E'_m(\R^d).$ 

For the singular GSP of Kondratiev-type $Q \in \mathcal{E}'(\R^d) \hat{\otimes}_\pi (S)_{-1}$, the coefficients $q_\gamma$, $\gamma \in \mathcal{I},$ in its chaos expansion form are elements of $\E_m'(\R^d)$ for some fixed $m \in \N_0$, and for the critical exponent  $p_Q$ we have
\begin{equation*}
	\normtri{Q}_{-p}^2 =\sum_{\gamma \in \mathcal{I}} \|q_{\gamma}\|^2_{\E'_m} (2\N)^{-p\gamma} < \infty, \quad p \geq p_Q.
\end{equation*}
\end{ex}
In this work we will use another  approach
based on Sobolev norms; see Section \ref{sec:singularQ_DefExi}.

\section{Stochastic parabolic equations with bounded potentials} \label{sec:boundedQ}
When 
considering the problem \eqref{Eq:ISH3}  with a potential $Q$, which is a GSP of Kondratiev-type, one faces  the problem of multiplication of GSPs.
%
%
To overcome the multiplication problem for random variables the Wick product is introduced  in \cite{HOUZ1996}, while the Wick product for GSPs  is introduced in \cite{LPSZ2015a}. 
Recall, for a normed vector space $X$,  and arbitrary GSPs $U,V \in X \hat{\otimes}_\pi (S)_{-p}$, $p\in \N$, with chaos expansions 
$U= \sum_{\alpha\in \mathcal{I} }u_\alpha   H_\gamma$ and $V= \sum_{\beta\in \mathcal{I}} v_\beta  H_\beta$ respectively,  the Wick product $U \lozenge V$ is defined by
$$U \lozenge V =\sum_{\gamma \in \mathcal I} \left( \sum_{\alpha + \beta = \gamma} u_\alpha  v_\beta  \right)  H_\gamma, \quad \omega \in \Omega. $$

For our analysis, for the operator $ \mathcal L $ we assume Assumption \ref{op_L} and define its action on  
a GSP  of Kondratiev-type $U\in X \hat{\otimes}_\pi (S)_{-1}$ with the chaos expansion  
	\begin{equation}
	\label{U chaos}
	U(t, x, \omega) = \sum\limits_{\gamma\in \mathcal I} u_\gamma(t, x) \, H_\gamma(\omega), 
\end{equation} 
by
$${\mathcal L} U(t, x, \omega) := \sum_{\gamma\in \mathcal I} \mathcal{L} \, u_\gamma(t, x) \, H_\gamma(\omega),$$
and we impose conditions on the given GSPs, the  force term $F$ and  the initial data $G$, which will together with the Assumption \ref{op_L}, remain in effect throughout the paper. 

\begin{ass} \label{as_F}
	The force term $F\in AC([0,T]; L^2(\R^d)) \hat{\otimes}_\pi (S)_{-1}$  is $AC([0,T]; L^2(\R^d)) $-valued GSP of Kondratiev-type  with the chaos expansion   
\begin{equation} \label{F chaos}
	F(t,x,\omega) = \sum_{\gamma \in \mathcal{I}} f_{\gamma}(t,x) H_{\gamma} (\omega),  \quad   f_\gamma \in AC([0,T];L^2 (\R^d)),
\end{equation} and with critical exponent $p_F \in \N$ such that
\begin{equation*} 
	\normtri{F}^2_{-p_F} = \sum\limits_{\gamma\in \mathcal I} \|f_\gamma\|^2_{AC([0,T];L^2 (\R^d))}  \,\, (2\mathbb N)^{-p_F\gamma}< \infty.
\end{equation*}
\end{ass}
\begin{ass} \label{as_G}
The initial condition $G\in D \hat{\otimes}_\pi (S)_{-1}$  is a $D$-valued GSP Kondratiev-type with the chaos expansion  
\begin{equation} \label{G chaos}
	G(x,\omega) = \sum_{\gamma \in \mathcal{I}} g_{\gamma}(x) H_{\gamma} (\omega),  \quad   g_\gamma \in D \subseteq L^2(\mathbb R^d),
\end{equation}
and with critical exponent $p_G \in \N$ such that  
\begin{equation*}
	\normtri{G}^2_{-p_G}=\sum\limits_{\gamma\in \mathcal I} \|g_\gamma\|_{L^2}^2  \, \, (2\mathbb N)^{-p_G\gamma}< \infty.
\end{equation*}
\end{ass}
As for the potential, in this section we assume boundedness in space of the GSP $Q$.
\begin{ass} \label{as_Qbnd}
The potential $Q \in L^\infty(\mathbb R^d) \hat{\otimes}_\pi (S)_{-1}$ is $L^\infty(\mathbb R^d)$-valued GSP of Kondratiev-type with the chaos expansion 
\begin{equation} \label{Q chaos}
Q(x, \omega) = \sum_{\gamma\in \mathcal I} q_\gamma( x) \, H_\gamma(\omega), \quad 
q_\gamma \in L^\infty(\mathbb R^d), 
\end{equation} 
and with critical exponent $p_Q \in \N$ such that
\begin{equation*} 
	\normtri{Q}^2_{-p_Q} = \sum_{\gamma \in \mathcal{I}}  \|q_\gamma\|^2_{L^\infty} (2\N)^{-p_Q \gamma} < \infty.
\end{equation*}
In addition, there exists $q>0$ such that
\begin{equation}\label{qLinf}
q :=  \sup_{\gamma \in \mathcal{I}}\|q_{\gamma}\|_{L^\infty} < \infty.
\end{equation}
\end{ass}
Under these assumptions we will be able to prove the existence of the stochastic weak solution in the following sence.
\begin{de}[Stochastic weak solution] \label{stweakSol}
  Generalized stochastic process $U$ of Kondratiev-type 
  is a \emph{stochastic weak solution} to the problem \eqref{Eq:ISH3}
  if $U \in AC([0,T], D) \hat{\otimes}_\pi (S)_{-1} $
    and the coefficients from its chaos expansion form \eqref{U chaos},  $u_\gamma \in AC([0,T], D),$ $\gamma \in \mathcal{I},$  satisfy deterministic system corresponding to the equation \eqref{Eq:ISH3}, i.e., $u_\gamma$ is a \em{weak solution} to the problem
  \begin{equation} \label{detProb}
  \left( \frac{\partial}{\partial t}  -  {\mathcal L} \right)
  u_\gamma (t,x) +q_{\mathbf{0}} (x) u_\gamma (t,x) = \tilde{f}_\gamma(t,x), \quad u_\gamma (0,x) = g_\gamma (x),
  \end{equation}
  where 
  $$\tilde{f}_\gamma(t,x) = \left\{ 
  \begin{array}{rl}
  	f_\gamma (t,x) -  \sum_{\substack{\alpha + \beta =\gamma \\ \alpha \not=\mathbf{0}}} q_\alpha (x) u_\beta (t,x), & \text{for } |\gamma|\geq 1, \\
  	f_\gamma(t,x), & \text{for } |\gamma|=0,
	\end{array}
	\right. $$
  and $f_\gamma$ and $g_\gamma$ are the coefficients in the chaos expansion forms \eqref{F chaos} and \eqref{G chaos} of generalized stochastic processes $F$ and $G$, respectively.
  \end{de}
 
 In \cite{GLO2023} the chaos expansion method is used to prove existence and uniqueness of solutions to \eqref{Eq:ISH3}, where also  estimates on the coefficients in the chaos expansion representation of the solutions are provided \cite[Theorem 4]{GLO2023}. 
 Here we state and prove a modified result. 
 
\begin{thm} 
	\label{T:genISH2} 
	Let $\mathcal L$, $Q$, $F$ and $G$ be as in Assumptions \ref{op_L}, \ref{as_Qbnd}, \ref{as_F}, \ref{as_G}, respectively.
	Let $\tilde{M}(T)$ be given by \eqref{mtilda}. 
	%
	%
	Then, 
	there exists a unique  generalized stochastic process  
	$$U \in AC([0,T];D)  \hat{\otimes}_\pi (S)_{-p_U} \subseteq AC([0,T];L^2 (\R^d)) \hat{\otimes}_\pi (S)_{-1} $$ 
	with chaos expansion \eqref{U chaos} 
	which is a stochastic weak solution  to
	the stochastic evolution initial value problem \eqref{Eq:ISH3} in sense of Definition \ref{stweakSol}.
For all $t\in [0,T]$ the coefficients $u_\gamma$, $\gamma \in \mathcal{I}$ satisfy
	\begin{equation}
	\label{EstKoefCE}
	\norm{u_\gamma (t,\cdot)}{L^2}  \leq M(t) \Bigg(a_\gamma(t) + \sum\limits_{k=1}^{|\gamma|} \big(\tilde{M}(t) q)^k \sum_{\substack{0\leq |\beta|\leq \gamma- k\\\beta < \gamma}} a_\beta(t) \Bigg),
	\end{equation}
	where $a_\gamma (t):= \|g_{\gamma} \|_{L^2 } + t \| f_{\gamma}\|$, $\gamma\in\mathcal I$,
	with $\|f_\gamma\|:= \sup_{t \in [0,T]} \|f_\gamma(t, \cdot)\|_{L^2}.$ 
	
	Moreover, for $p\geq p_U:=  \big\lfloor \max\left\{2mp_F,2mp_G,\frac{m(3+s)}{m-1}\right\}\big\rfloor +1 $, $m \in \N\setminus\{1\}$ arbitrary, where $s = s(M, w,T, q)$ is given by  $s = \frac{\ln{(\tilde{M}(T) q)^{2}}}{\ln 2} +1$ for $(\tilde{M}(T) q)^2>1$ and $s=0$ for $(\tilde{M}(T) q)^2\leq 1$ , it holds
	\begin{equation} \label{ocena}
	\normtri{U}^2_{-p} \leq 3 M(T)^2 A\left( 1 + 2C_{\frac{p}{2m}} C_{\frac{p(m-1)}{m}-s-2} \right),
	\end{equation}
	where
	\begin{equation} \label{A}
A:= \normtri{G}_{-p_G}^2 + T^2 \normtri{F}_{-p_F}^2.
\end{equation}
\end{thm}
%
\begin{rmk}
	Comparing  conditions of Theorem 4 in  \cite{GLO2023},
	we remark  that the conditions were stronger:
	the condition on the coefficients in the chaos expansion of the potential $Q$  was given by 
%
	$\|q_{\gamma}\|_{L^\infty} \leq \|q_\mathbf{0}\|_{L^\infty}$
	for all $\gamma \in \mathcal{I}$, 
while the condition  $\tilde{M}(T) \|q_\mathbf{0}\|_{L^\infty} \not= 1$ is  now dropped.
\end{rmk}


\proof
As in proof of Theorem 4 in \cite{GLO2023}, given stochastic processes $Q$, $F$ and $G$ appearing in the problem \eqref{Eq:ISH3} are represented in the chaos expansion forms \eqref{Q chaos}, \eqref{F chaos}, \eqref{G chaos}, respectively, and assuming the solution $U$ in the chaos expansion form
\eqref{U chaos}, by the definition of the Wick product and the chaos expansion method  the initial problem \eqref{Eq:ISH3} is reduced to a system of deterministic parabolic equations of the form \eqref{detProb}. By Theorem \ref{T:DetEq} there exists 
a unique solution $u_{\gamma} \in AC([0,T],D)$  to the corresponding deterministic problem, and it satisfies 
\begin{eqnarray}\label{estFromThm}
	\norm{u_\gamma (t,\cdot)}{L^2}  & \leq  &
	M(t)\left( a_\gamma (t) +  \sum\limits_{\mathbf{0} \leq \beta < \gamma} \linf{q_{\gamma-\beta}} \int_0^t \ltwo{ u_{\beta} (s, \cdot)}\,ds\right) \nonumber \\
	 & \leq  &	M(t)\left( a_\gamma (t) +  \sum\limits_{\mathbf{0} \leq \beta < \gamma} q \int_0^t \ltwo{ u_{\beta} (s, \cdot)}\,ds\right),
\end{eqnarray}
where $q$ is defined by \eqref{qLinf}.
By induction,  one proves that coefficients $u_\gamma$, $\gamma \in \mathcal{I}$, satisfy the estimate \eqref{EstKoefCE}. 
Indeed, since the estimate for $|\gamma|=0$
is clearly true, then  assuming the estimate \eqref{EstKoefCE} for every $\beta \in \mathcal{I}$ 
with $|\beta|\leq n$, namely
\begin{equation} \label{ih}
		\norm{u_\beta (t,\cdot)}{L^2}  \leq M(t) \Bigg(a_\beta(t) + \sum\limits_{k=1}^{|\beta|} \big(\tilde{M}(t) q)^k \sum_{\substack{0\leq |\alpha|\leq \beta- k\\\alpha < \beta}} a_\alpha(t) \Bigg),
\end{equation}
one  shows that \eqref{EstKoefCE} holds for $\gamma \in \mathcal{I}$ with $|\gamma|=n+1$  as follows.
First, integrating  \eqref{ih}, and using Lemma \ref{Lema M(t)}, one obtains
$$
\int_0^t \ltwo{ u_{\beta} (s, \cdot)}\,ds \leq a_\beta (t) \tilde{M}(t) + \sum\limits_{l=1}^{|\beta|} \tilde{M}(t)^{l+1} q^l  \sum_{\substack{0 \leq |\alpha| \leq |\beta|-l \\ \alpha<\beta}} a_\alpha (t),
$$ 
and then, starting from  \eqref{estFromThm}  that 
\begin{align*}
	\norm{u_\gamma (t,\cdot)}{L^2}& \leq 
	M(t) \Bigg\{  a_\gamma (t)    
	+  \sum\limits_{\mathbf{0} \leq \beta < \gamma} q\, a_\beta (t) \tilde{M}(t) 
	+  \sum\limits_{\mathbf{0} < \beta < \gamma} q \sum\limits_{l=1}^{|\beta|}\tilde{M}(t)^{l+1} q^l  			\sum_{\substack{0 \leq |\alpha| \leq |\beta|-l \\ \alpha<\beta}} a_\alpha (t)   \Bigg\}\\
	& \leq M(t) \Bigg\{  a_\gamma (t)    
	+  \sum\limits_{\mathbf{0} \leq \beta < \gamma} a_\beta (t) \tilde{M}(t) q
	+ \sum_{l=1}^n (\tilde{M}(t)q)^{l+1} \sum\limits_{\mathbf{0} < \beta < \gamma} \Bigg( \sum_{\substack{\mathbf{0} \leq \alpha<\beta\\ 0 \leq |\alpha| \leq |\beta|-l }} a_\alpha (t)\Bigg)\Bigg\}.
\end{align*}
Next, merging the sums over $\mathbf{0} < \beta < \gamma$ (which implies $|\beta| \leq |\gamma|-1$) and over $\mathbf{0} \leq  \alpha<\beta$ with $ |\alpha| \leq |\beta|-l $ into the sum over $\mathbf{0} \leq \alpha <  \gamma$ with $ |\alpha| \leq |\gamma|-(l+1)$ one obtains
$$
	\sum_{l=1}^n (\tilde{M}(t)q)^{l+1}\!\!\! \sum_{\substack{\mathbf{0} \leq \alpha <\gamma \\ 0 \leq |\alpha| \leq |\gamma|-(l+1) }} \!\!\! a_\alpha (t) 
	= \sum_{k=2}^{|\gamma|} (\tilde{M}(t)q)^{k} \!\!\! \sum_{\substack{  0 \leq |\alpha| \leq |\gamma|-k \\\alpha <\gamma }} a_\alpha (t),
$$
and therefore
\eqref{EstKoefCE} follows.\\

To show that the solution $U$ is a GSP of Kondratiev-type, one has to show that there exists a critical exponent $p_U>1$ (to be determined bellow) such that the sum $$\normtri{U}_{-p}^2  := \sum_{ \gamma\in \mathcal{I}}\|u_\gamma\|^2_{AC([0,T]; L^2(\R^d))} (2\N)^{-p\gamma}$$  is finite for all $p \geq p_U$. 
Using  \eqref{EstKoefCE} and the inequality 
\begin{equation}\label{kvadratsume}
	\left(\sum_{k=1}^{n} x_k \right)^2 \leq n \sum_{k=1}^{n}x_k^2, \quad x_k \in \R,
\end{equation}
 we find that $\normtri{U}^2_{-p} $ is bounded by
$$3 M(T)^2 \Bigg( \sum_{\gamma \in \mathcal{I}}\| g_{\gamma}\|_{L^2}^2 (2 \N)^{-p \gamma}  + T^2 \sum_{\gamma \in \mathcal{I}} \| f_{\gamma}\|^2 (2 \N)^{-p \gamma} + \sum_{\gamma \in \mathcal{I}} \Big(\sum\limits_{k=1}^{|\gamma|} (\tilde{M}(T)q)^k  \sum_{\substack{0 \leq |\beta| \leq |\gamma|-k \\ \beta<\gamma}} a_\beta (T) \Big)^2 (2 \N)^{-p \gamma} \Bigg).$$ 
Choosing $p \geq \max\{p_F,p_G\}$ 
one estimates sum of first two terms by a finite term $A$ defined by \eqref{A}.
To estimate the third term 
$$S_3:=\sum_{\gamma \in \mathcal I} \Big(\sum\limits_{k=1}^{|\gamma|} \big(\tilde{M}(T) q)^k \sum_{\substack{0\leq |\beta|\leq \gamma- k\\\beta < \gamma}} a_\beta(T) \Big)^2 (2\N)^{-p\gamma}$$ 
 we first apply the inequality \eqref{kvadratsume}  and use $(2\N)^{-p\gamma} = (2\N)^{-\frac{p}{m}\gamma}\cdot (2\N)^{-\frac{p(m-1)}{m}\gamma}$, for arbitrary chosen $m \in \N$, $m>1$,  and obtain
\begin{align}
	S_3 &\leq \sum_{\gamma \in \mathcal I} |\gamma| \sum\limits_{k=1}^{|\gamma|}
		(\tilde{M}(T) q)^{2k} \big( \sum_{\substack{0\leq |\beta|\leq \gamma- k\\\beta < \gamma}} a_\beta(T)  \, (2\N)^{-\frac{p}{2m}\gamma}\big)^2 (2\N)^{-\frac{p(m-1)}{m}\gamma}\nonumber \\
		&\leq  \sum_{\gamma \in \mathcal I} |\gamma| \sum\limits_{k=1}^{|\gamma|}
		(\tilde{M}(T) q)^{2k} \big( \sum_{\beta \in \mathcal I} a_\beta(T)  \, (2\N)^{-\frac{p}{2m}\beta}\big)^2 (2\N)^{-\frac{p(m-1)}{m}\gamma}.\label{S3 pocetak}
\end{align}
The inner sum we estimate  using the Cauchy-Schwarz inequality \eqref{CSmulti}

\begin{equation*}
	\begin{split}\Big(\sum_{\beta \in \mathcal I} a_\beta(T)  \, (2\N)^{-\frac{p}{2m}\beta} \Big)^2&= \Big(\sum_{\beta \in \mathcal I} \big(a_\beta(T)  \, (2\N)^{-\frac{p}{4m}\beta} \big) (2\N)^{-\frac{p}{4m}\beta} \Big)^2 \\&\leq 
		\Big(\sum_{\beta \in \mathcal I} a_\beta^2(T)  \, (2\N)^{-\frac{p}{2m}\beta} \Big) \cdot \Big(\sum_{\beta \in \mathcal I} (2\N)^{-\frac{p}{2m}\beta}\Big)\\
		& = C_{\frac{p}{2m}} \cdot  \sum_{\beta \in \mathcal I} a_\beta^2(T)  \, (2\N)^{-\frac{p}{2m}\beta} 
	\end{split}
\end{equation*} 
with the constant  $C_{\frac{p}{2m}}$, defined  via \eqref{Cp}, which is finite for $p>2m$ by Lemma \ref{L:konvCp}.  
For $p\geq \max\{2m p_F, 2m p_G\}$, we have 
\begin{align*}
	\sum_{\beta \in \mathcal I} a^2_{\beta} (T)(2\N)^{-\frac{p\beta}{2m}} 
	& =\sum_{\beta \in \mathcal I} (\|g_{\beta}\|_{L^2} + T  \| f_{\beta}\|)^2 (2\N)^{-\frac{p\beta}{2m}} \leq 2 \, \sum_{\beta \in \mathcal I}  \left(  \|g_{\beta}\|_{L^2}^2 +   T^2 \|f_{\beta}\|^2\right) (2\N)^{-\frac{p\beta}{2m}}\nonumber\\
	& = 2 \sum_{\beta \in \mathcal I}   \|g_{\beta}\|_{L^2}^2 (2\N)^{-\frac{p\beta}{2m}} + 2  T^2 \sum_{\beta \in \mathcal I}    \|f_{\beta}\|^2 (2\N)^{-\frac{p\beta}{2m}} \nonumber\\
	&\leq
	2 \sum_{\beta \in \mathcal I}   \|g_{\beta}\|_{L^2}^2 (2\N)^{-p_G\beta} + 2  T^2 \sum_{\beta \in \mathcal I}    \|f_{\beta}\|^2 (2\N)^{-p_F\beta}= 2A, \nonumber\\	
\end{align*}
where  
 $A$ is given by \eqref{A}, and it is finite by  convergence conditions for processes $F$ and $G$. 
Thus, for $p\geq \max\{2 m p_F, 2 m p_G\}$,
\begin{equation}
	\Big(\sum_{\beta \in \mathcal I} a_\beta(T)  \, (2\N)^{-\frac{p}{2m}\beta} \Big)^2 \leq 2AC_{\frac{p}{2m}},
\end{equation} 
which, used in \eqref{S3 pocetak}, gives final estimate for the third term 
\begin{equation*}
		S_3\leq 2AC_{\frac{p}{2m}} \sum_{\gamma \in \mathcal I} |\gamma| \sum\limits_{k=1}^{|\gamma|}
		(\tilde{M}(T) q)^{2k}  (2\N)^{-\frac{p (m-1)}{m}\gamma}.
\end{equation*}
To see that $\normtri{U}^2_{-p}$ is finite, we note that the term $\sum\limits_{k=1}^{|\gamma|}
(\tilde{M}(T) q)^{2k} $ is the sum of finitely many elements (exactly $|\gamma|$ of them) of a geometric sequence with the first term $r=(\tilde{M}(T) q)^2$. 

If $r\leq1$, the sum $\sum\limits_{k=1}^{|\gamma|}
(\tilde{M}(T) q)^{2k} = \sum\limits_{k=1}^{|\gamma|} r^k $ is  estimated as $\sum_{k=1}^{|\gamma|} r^k \leq |\gamma|$, and by Lemma \ref{lemanew} with $d(\gamma)=|\gamma|^2$,  $S_3 $ could be estimate as
\begin{equation*}
	\begin{split}
		S_3  &\leq 2AC_{\frac{p}{2m}} \sum_{\gamma \in \mathcal I} |\gamma| \sum\limits_{k=1}^{|\gamma|}
		r^k\,   (2\N)^{-\frac{p (m-1)}{m}\gamma} \leq  2AC_{\frac{p}{2m}} \sum_{\gamma \in \mathcal I} |\gamma|^2  (2\N)^{-\frac{p (m-1)}{m}\gamma}\\
		& =  2 AC_{\frac{p}{2m}} D_{\frac{p (m-1)}{m}}
		\leq 
		2AC_{\frac{p}{2m}}C_{\frac{p (m-1)}{m} -2 }, 
	\end{split}
\end{equation*}
which is finite for $p > \max\{2mp_F,2mp_G,\frac{3m}{m-1}\}$ by Lemma \ref{L:konvCp}. 

If $r>1$ then we estimate each term in the sum by the largest term and get
%
$	\sum_{k=1}^{|\gamma|} r^k \leq  |\gamma|  r^{|\gamma|}$.
 By Lemma \ref{lemanew} with $d(\gamma) = |\gamma|^2 r^{|\gamma|}$, there exists $s>0$ such that
\begin{equation*}
		S_3 
		 \leq  2AC_{\frac{p}{2m}} \sum_{\gamma \in \mathcal I} |\gamma|^2 r^{|\gamma|}  (2\N)^{-\frac{p(m-1)}{m}\gamma} =  2AC_{\frac{p}{2m}} D_{\frac{p(m-1)}{m}} \\
		\leq 
		 2AC_{\frac{p}{2m}}C_{\frac{p(m-1)}{m}-s-2},
\end{equation*}
which is finite for $p>\max\left\{2mp_F,2mp_G,\frac{m(3+s)}{m-1}\right\}.$

Therefore, $\normtri{U}^2_{-p}$ is finite for  $p>\max\left\{2mp_F,2mp_G,\frac{m(3+s)}{m-1}\right\}$, where $s = 0$
in the case $r\leq 1$,  and otherwise  $s \geq 0$ is chosen such that $r^{|\gamma|}\leq (2\N)^{s\gamma}$.
We note that $s$ is a constant depending on $r$, i.e., on $\tilde{M} (T)$ and $q$, i.e., on $M$, $w$, $T$ and $q$. Thus, critical exponent in chaos expansion for the solution $U$ is $p_U >\max\left\{2mp_F,2mp_G,\frac{m(3+s)}{m-1}\right\}$, with $s = s(M, w, T, q)$.
More precisely, 
$$\normtri{U}^2_{-p} < 3 M(T)^2 A\left( 1 + 2C_{\frac{p}{2m}} C_{\frac{p(m-1)}{m}-s-2}\right),$$ 
is finite for $p\geq p_U$ with  
$$ p_U:=  \big\lfloor \max\left\{2mp_F,2mp_G,\frac{m(3+s)}{m-1}\right\}\big\rfloor +1 $$
and where $s=0$ for $r \leq 1$ and $s := \frac{\ln{(\tilde{M}(T) q)^{2}}}{\ln 2} +1$,  for $r>1$, by Remark \ref{Rem_on_s}.
%
%

The uniqueness of the solution $U$ follows from the uniqueness of its coefficients $u_\gamma$, $\gamma\in \mathcal I,$ and the uniqueness of  the chaos expansion representation in terms of the Fourier--Hermite stochastic polynomials  that form an orthogonal basis of the space $L^2(\mu)$.\endproof

\section{Stochastic parabolic equations with singular potentials} \label{sec:singularQ_DefExi}

Finally, we are ready to consider the stochastic parabolic equations of the form \eqref{Eq:ISH3} with singular space dependent random potentials $Q\in H^{-s}_c(\R^d) \hat{\otimes}_\pi (S)_{-1}$, $s>0$, where $ H^{-s}_c(\R^d)$ is Sobolev space of compactly supported distributions defined in the Introduction. 
Namely, we consider the problem 
\begin{equation*}
\left(\frac{\partial}{\partial t}  - \mathcal L \right) \, U(t, x, \omega) + Q(  x, \omega) \, \lozenge \, U(t, x, \omega) = F(t, x, \omega),\quad
U(0, x, \omega) = G(x, \omega),
\end{equation*} for $t\in (0,T],$ $x\in \mathbb R^d$, $\omega\in \Omega$, supposing Assumption \ref{op_L} for the operator $\mathcal L$, and  Assumptions \ref{as_F} and \ref{as_G} for  the load vector $F$ and the initial data $G$.  For the potential $Q$, the Assumption \ref{as_Qbnd} from Section \ref{sec:boundedQ} will be replaced with the following one.
\begin{ass} \label{as_Qsing}
The given potential is singular GSP of Kondratiev-type, i.e., $Q \in H^{-s}_c(\R^d) \hat{\otimes}_\pi (S)_{-1}$, $s>0$, 
	with the chaos expansion representation
	\begin{equation}
		\label{q chaos}
	Q(x,\omega) =
	 \sum\limits_{\gamma\in \mathcal I } q_\gamma (x)  \, H_\gamma (\omega), \quad q_\gamma \in  H^{-s}_{c} (\mathbb R^d),
	\end{equation}
such that for critical exponent $p_Q\in \N$ the following growth condition holds
	$$
	\normtri{Q}^2_{-p_Q}= \sum_{\gamma \in \mathcal I} \| q_{\gamma}\|_{H^{-s}_{c}}^2 (2\N)^{-p_Q \gamma}  
 < \infty,
  $$	 
  where  
$\| \cdot\|_{H^{-s}_{c}}$ is seminorm given by \eqref{H}.
 In addition, we suppose that there exists  $q>0$ such that  
	\begin{equation}\label{noviuslov}
		q :=  \sup_{\gamma \in \mathcal{I}}\|q_{\gamma}\|_{H^{-s}_c}<\infty.
	\end{equation}
\end{ass}
%
%

\subsection{Regularization of the singular GSPs of Kondratiev-type}
When dealing with stochastic potentials $Q$ that are singular with respect to 
$x$, one encounters the fundamental issue of multiplying distributions, which is not well-defined in the classical sense. To address this problem, a common approach involves regularizing the singular GSP 
$Q$ with respect to 
$x$ using a classical regularization procedure.
In this work, we tackle this challenge by employing a regularization technique based on convolution with a mollifier. The mollifier smooths out the singularities of 
$Q$, thereby enhancing its mathematical properties and enabling a rigorous analysis of the problem \eqref{Eq:ISH3}. This approach serves as a foundational step toward studying stochastic parabolic equations with singular potentials in a mathematically consistent framework.


Let  $\varphi$ be a \emph{mollifier}, i.e., a compactly supported smooth function  $\varphi \in \D (\mathbb{R}^d)$, with  $\varphi \geq 0$ and $\int_{\mathbb{R}^d} \varphi (x)\,dx=1$, and let  $(\varphi_{\varepsilon})_{\varepsilon\in (0,1]}$  be the \emph{mollifying net} defined as
\begin{equation}\label{molinet}
\varphi_{\varepsilon} (x):= \dfrac{1}{\varepsilon^d} \, \varphi \left( \dfrac{x}{\varepsilon}\right) \in \D(\mathbb{R}^d).
\end{equation}
For a given singular stochastic potential $Q\in\E'(\R^d) \hat{\otimes}_\pi (S)_{-1}$ with the chaos expansion \eqref{q chaos}, a
regularization with respect to $x$ via convolution is the process in which we convolve the potential $Q$, with respect to $x$,  and the mollifying  net  $(\varphi_\varepsilon)_{\varepsilon \in (0,1]}$ given by \eqref{molinet}, and obtain the \emph{regularizing net}  of generalized stochastic processes
with compactly supported  smooth 
 coefficients, i.e., net  
$(Q_\varepsilon)_{\varepsilon\in (0,1]}\in \D(\R^d) \hat{\otimes}_\pi (S)_{-1} $. More precisely, for each $\varepsilon > 0$ the process $Q_\varepsilon$  is defined as
\begin{equation}\label{reg_Q}
Q_\varepsilon:= Q\ast_x \varphi_\varepsilon=
\sum\limits_{\gamma\in \mathcal{I}} (q_\gamma\ast_x\varphi_\varepsilon) \, H_\gamma,  
\end{equation}
 where $q_\gamma \in \E'(\R^d)$ implies $q_\gamma \ast_x\varphi_\varepsilon\in \D(\mathbb R^d)$  for all $\gamma \in \mathcal{I}$.

In the sequel, the key concept is {\em moderateness}, see \cite{GKOS2001,MO1992}. This concept is particularly natural in the context of regularizations, as regularizations of distributions inherently produce moderate families. Specifically, structure theorems for distributions establish that compactly supported distributions 
$\mathcal{E}'(\R^d)$ can be regarded as a subset of 
$C^\infty (\R^d)$-moderate families.

Building on this foundation, we proceed to define moderate nets of generalized stochastic processes of Kondratiev-type, which will serve as a crucial tool for analyzing singular stochastic potentials and their regularizations.
\begin{de}[Moderate net of GSPs]\label{ModGSp2}
	Let $({X},\| \cdot\|_{X})$ be a Banach (or Frech\'et) space and let  $(U_\varepsilon)_{\varepsilon \in (0,1]}\subset {X} \hat{\otimes}_\pi (S)_{-1}$ be a  net of GSPs of Kontratiev-type.
	%
	The net  
	$(U_\varepsilon)_{\varepsilon \in (0,1]}$ is ${X} \hat{\otimes}_\pi (S)_{-1}$\emph{-moderate}
	if for each $\varepsilon>0$, there exist $N \in \mathbb{N}_0$   and $C>0$, independent of $\varepsilon$, and there exists $p_{\varepsilon} \in \N$, so that 
	$$\normtri{U_{\varepsilon}}_{-{p}}
	\leq C \varepsilon^{-N}, \quad p\geq p_{\varepsilon}.$$  
\end{de}

\begin{lem} \label{L:Qeps}
	Let $Q \in \E'(\mathbb{R}^d) \hat{\otimes}_\pi (S)_{-1}$ with the chaos expansion \eqref{q chaos}
	and  let $(Q_\varepsilon)_{\varepsilon\in (0,1]}$ be a regularizing net defined in \eqref{reg_Q}. Then 	$(Q_\varepsilon)_{\varepsilon\in (0,1]}$ converges to $Q$ in $\E'(\R^d) \hat{\otimes}_\pi (S)_{-1}$ as $\varepsilon \to  0$. If $Q \in H^{-s}(\mathbb{R}^d) \hat{\otimes}_\pi (S)_{-1}$, for $s>0$, then
$(Q_\varepsilon)_{\varepsilon\in (0,1]}$ is $L^\infty(\mathbb R^d) \hat{\otimes}_\pi (S)_{-1}$-moderate. 

\end{lem}
\proof Let $\phi \in \mathcal{E}(\R^d)$ be an arbitrary test function. Let us show that $\langle Q_\eps (x,\omega), \phi(x) \rangle \to \langle Q(x,\omega), \phi (x)\rangle$ as $\eps \to 0.$ 
It is known that $\langle q_\gamma \ast_x \varphi_\eps, \phi \rangle \to \langle q_\gamma, \phi\rangle$ as $\eps \to 0.$ 
%
%
%
By Lebesgue dominated convergence theorem, with a dominant function $g(x):=\sup_{\eps \in (0,1]}\sup_{\gamma \in \mathcal{I}} |\langle q_\gamma,\varphi_\eps \ast_x \phi(x)\rangle|,$
we have
\begin{eqnarray*}
	\lim_{\eps\to 0} \langle Q_\varepsilon(x,\omega), \phi(x) \rangle &= &\lim_{\eps\to 0} 	\sum\limits_{\gamma\in \mathcal{I}}  \langle q_\gamma\ast_x\varphi_\varepsilon( x), \phi(x) \rangle \, H_\gamma(\omega) \\
	& =& 	\sum\limits_{\gamma\in \mathcal{I}} \lim_{\eps\to 0} \langle (q_\gamma\ast_x\varphi_\varepsilon)( x), \phi(x) \rangle \, H_\gamma(\omega)\\
	&=& \sum\limits_{\gamma\in \mathcal I} \langle q_\gamma( x),\phi(x) \rangle \, H_\gamma(\omega) \\
	&=& \langle Q(x, \omega), \phi(x) \rangle.
\end{eqnarray*} 
Since for every $\gamma \in \mathcal{I}$, $q_\gamma \in H^{-s}_{c}(\R^d)$, for fixed $s>0$,  we have  that $q_\gamma = \sum_{|\alpha| \leq s} \left(\frac{\partial}{\partial x}\right)^\alpha k_{\gamma,\alpha}$,  $k_{\gamma,\alpha} \in L^2(\R^d),$ and 
\begin{equation*} 
q_\gamma \ast \varphi_\eps = \sum_{|\alpha| \leq s} \left(\frac{\partial}{\partial x}\right)^\alpha k_{\gamma,\alpha} \ast \varphi_\eps = \sum_{|\alpha| \leq s}  k_{\gamma,\alpha} \ast  \left(\frac{\partial}{\partial x}\right)^\alpha \varphi_\eps = \sum_{|\alpha| \leq s} \dfrac{1}{\eps^{|\alpha|}} k_{\gamma,\alpha} \ast (\partial^\alpha \varphi)_\eps.
\end{equation*}
By Young's inequality (see \cite[Proposition 8.9, p. 241]{Folland1999}) 
\begin{equation*} 
	\|k_{\gamma,\alpha} \ast (\partial^\alpha \varphi)_\eps\|_{L^\infty} \leq \|k_{\gamma,\alpha}\|_{L^2} \|(\partial^\alpha \varphi)_\eps\|_{L^2} = \dfrac{1}{\varepsilon^{\frac{d}{2}}} \|k_{\gamma,\alpha}\|_{L^2} \|\partial^\alpha \varphi \|_{L^2},
	\end{equation*}
and so
\begin{eqnarray*} 
	\|q_\gamma \ast \varphi_\eps\|_{L^\infty} & \leq & \sum_{|\alpha| \leq s} \dfrac{1}{\eps^{|\alpha| + \frac{d}{2}}} \|k_{\gamma,\alpha}\|_{L^2} \|\partial^\alpha \varphi \|_{L^2}  \leq  \dfrac{1}{\eps^{s + \frac{d}{2}}} C_\varphi \sum_{|\alpha| \leq s}  \|k_{\gamma,\alpha}\|_{L^2} \\
	&\leq & \dfrac{C_\varphi \lceil s\rceil^{\frac{1}{2}}}{\eps^{s + \frac{d}{2}}} \left( \sum_{|\alpha| \leq s} \|k_{\gamma,\alpha}\|_{L^2} \right)^{\frac{1}{2}}.
\end{eqnarray*}
This estimate is true for any choice of $k_{\gamma,\alpha}$, and thus also for the one reading the infimum, i.e.,
\begin{equation}\label{star1}
		\|q_\gamma \ast \varphi_\eps\|_{L^\infty}  \leq   \dfrac{C_\varphi \lceil s\rceil^{\frac{1}{2}}}{\eps^{s + \frac{d}{2}}} \|q_\gamma\|_{H^{-s}_c}, 
\end{equation}
yielding 
\begin{equation} \label{normQeps}
\normtri{Q_{\eps}}^2_{p_Q} = \sum_{\gamma \in \mathcal{I}} \| q_\gamma \ast \varphi_\eps\|^2_{L^\infty} (2\N)^{-p_{Q} \gamma} \leq  \dfrac{C_\varphi^2 \lceil s\rceil}{\eps^{s + \frac{d}{2}}} \sum_{\gamma \in \mathcal{I}}\|q_\gamma\|_{H^{-s}_c}^2 (2\N)^{-p_{Q}\gamma}  =  \dfrac{C_\varphi^2 \lceil s\rceil}{\eps^{s + \frac{d}{2}}} \normtri{Q}^2_{p_Q}.
\end{equation}
 Thus, with $N=s+\frac{d}{2}>0$, and $C=C_\varphi^2 \lceil s\rceil  \normtri{Q}^2_{p_Q}$, we have that 
$$\normtri{Q_\eps}_{-p_{Q}}^2 \leq C \eps^{-N},$$
i.e.,  $(Q_\varepsilon)_{\varepsilon\in (0,1]}$ is $L^\infty(\mathbb R^d) \hat{\otimes}_\pi (S)_{-1}$-moderate. 
\endproof
\begin{rmk} \label{Rem_on_reg}
\item[a)] From \eqref{normQeps} we also read that for each $\eps>0$ the critical exponent of the GSP $Q_\eps$ defined  by \eqref{reg_Q} is equal to the critical exponent of the singular GSP $Q$, i.e., $p_{Q_\eps} = p_{Q}$.
\item[b)] For $Q \in \mathcal{E}'(\R^d)$, we need additional assumptions in order to prove $L^\infty (\R^d)$-moderatness of net $(Q_\eps)_\eps$.
Namely,	by \cite[Lemma 4 (a)]{GLO2021} it follows that for each $\gamma \in \mathcal I$ the net 
	$(q_\gamma \ast_x \varphi_\eps)_{\varepsilon \in (0,1]}$ is $L^\infty(\mathbb R^d)$-moderate, i.e., there exist $N_\gamma \in \N$ and $C_\gamma>0$ such that $\|q_\gamma \ast_x \varphi_\eps\|_{L^\infty} \leq C_\gamma \eps^{-N_\gamma}.$ Therefore, we have
	$$\normtri{Q_\eps}_{-p_{Q}}^2 =\sum_{\gamma \in \mathcal{I}} \|q_\gamma \ast \varphi_\eps\|_{L^\infty}^2 (2\N)^{-p_{Q} \gamma} \leq \sum_{\gamma \in \mathcal{I}}  C_\gamma^2 \eps^{-2N_\gamma} (2\N)^{-p_{Q} \gamma}.$$
Therefore, we need to suppose the existence of $N_1:=\inf_{\gamma \in \mathcal{I}} N_\gamma <N_\gamma$ and $C_\gamma < C_1 := \sup_{\gamma \in \mathcal{I}} C_\gamma$, and then it follows
	$$\normtri{Q}_{-p_{Q}}^2 < C_1^2 \eps^{-2N_1}  \sum_{\gamma \in \mathcal{I}} (2\N)^{-p_{Q} \gamma} = C_1^2 \eps^{-2N_1}  C_{p_{Q}}.$$ Finally,  $\normtri{Q_\eps}_{-p_{Q}}^2 \leq C \eps^{-N} ,$ with $N=2N_1$ and $C=C_1^2 C_{P_Q}.$
	 
	Similarly, by \cite[Lemma 4 (a)]{GLO2021},  for each $\gamma \in \mathcal{I}$ the net $(q_\gamma \ast_x \varphi_\eps)_{\varepsilon \in (0,1]}$ is log-type $L^\infty(\mathbb R^d)$-moderate, i.e., there exists $C_\gamma > 0$ such that
	$\norm{q_\gamma \ast_x \varphi_\eps}{L^\infty}\leq C_\gamma \log{\frac{1}{\eps}},$ and therefore if we suppose that there exists $C_2=\sup_{\gamma \in \mathcal{I}} C_\gamma > C_\gamma$ then
	$$\normtri{Q_\eps}_{-p_{Q}}^2 
	\leq \log^2 \dfrac{1}{\eps} \, \sum_{\gamma \in \mathcal{I}} C_\gamma^2  (2\N)^{-p_{Q} \gamma} < C_2^2   C_{p_{Q}} \log^2 \dfrac{1}{\eps} < C_2^2 \delta \log^2 \dfrac{1}{\eps}.$$ 
	Thus,
	$$\normtri{Q_\eps}_{-p_{Q}} \leq C \log \dfrac{1}{\eps},$$ with $C = C_2\sqrt{\delta},$
	i.e., $(Q_\eps)_\eps$ is log-type moderate. 
\end{rmk}

\subsection{Stochastic very weak solution }\label{sec:4}
The concept of stochastic very weak solutions for parabolic SPDEs with deterministic singular potentials is introduced in \cite{GLO2021}. Two distinct approaches are proposed for defining such solutions: Solution concept 1 and Solution concept 2.
In the Solution concept 1, the chaos expansion method is applied directly to the original stochastic problem, yielding a system of deterministic equations with a singular potential. The coefficients of the chaos expansion are then obtained by solving these deterministic equations in the sense of very weak solutions.
In the Solution concept 2, the singular potential is first regularized using an appropriate regularization technique. The chaos expansion method is then applied to the resulting net of stochastic parabolic equations, leading to the formulation of the stochastic very weak solution.

In this work, we introduce the notion of stochastic very weak solutions for parabolic SPDEs with random potentials that are singular in space. To construct such solutions, we adopt Solution concept 2, which we recall below.
 \begin{de}[Stochastic very weak solution]\label{option 2} Let the potential $Q$ be in  {$\mathcal E'(\mathbb R^d) \hat{\otimes}_\pi (S)_{-1}$}.	A net  of stochastic processes $(U_{\varepsilon})_{\varepsilon\in (0,1]}$ in  $AC([0,T];L^2 (\R^d)) \hat{\otimes}_\pi (S)_{-1}$
	is  a {\rm very weak solution} to  the stochastic  problem 
	\eqref{Eq:ISH3}
	if  there exists an $L^\infty(\mathbb R^d)$-moderate  regularizing net of stochastic processes $(Q_\varepsilon)_{\varepsilon\in (0,1]}$ of the stochastic process $Q$ (which is a distribution with respect to $x$), such that for every $\varepsilon\in (0,1]$  the process $U_\varepsilon$ is a solution to the regularized stochastic parabolic  equation
	\begin{equation}
	\label{SP reg}
	\begin{split}\small 
	\left( \frac{\partial}{\partial t}  - \mathcal{L} \right) \, U(t, x,\omega) + {Q_\varepsilon (x,\omega)} \lozenge  U(t, x,\omega) &= F(t, x,\omega),	\\
	U(0, x,\omega) &= G(x,\omega),
	\end{split}
	\end{equation}
	and $(U_{\varepsilon})_{\varepsilon\in (0,1]}$  is a moderate net of GSPs of Kondratiev-type {in the sence of Definition \ref{ModGSp2}. }
\end{de}
\begin{thm}[{Existence of a stochastic very weak solution}] \label{T:main}
	Let the operator $\mathcal L$,  the  potential $Q$,
	the force term $F$, and the initial condition $G$ be as in Assumptions \ref{op_L}, \ref{as_Qsing}, \ref{as_F} and \ref{as_G}, respectively.
	Then,  the stochastic parabolic problem \eqref{Eq:ISH3}	 
	has a  stochastic very weak solution  $(U_\varepsilon)_{\varepsilon\in (0,1]}$  in  $AC([0,T];L^2(\R^d)) \hat{\otimes}_\pi (S)_{-1}$.
\end{thm}

\proof 
According to Lemma \ref{L:Qeps} by \eqref{reg_Q} is given an $L^\infty(\mathbb R^d)$-moderate  regularization $\net{Q_\eps}$ of the singular {GSP} $Q$.
The regularized problem for stochastic parabolic equation \eqref{Eq:ISH3} 
is a net of problems of the form {\eqref{SP reg}}.
 For each $\eps\in(0,1]$ the problem \eqref{SP reg} is of type \eqref{Eq:ISH3} with bounded potential GSP $Q_\varepsilon$  with bounded coefficients in chaos expansion representation $q_\gamma\ast \vphi_\eps \in L^\infty(\mathbb R^d)$. 
%
By Assumption \ref{as_Qsing} $q= \sup_{\gamma \in \mathcal{I}} \|q_\gamma\|_{H^{-s}_c}$ is finite, which implies that for each $\eps>0$
	$$q_\eps: = \sup_{\gamma \in \mathcal{I}}\|q_{\gamma}\ast \vphi_\eps\|_{L^\infty} <\infty.$$
Therefore, for each $\eps\in(0,1]$,  one can apply  Theorem \ref{T:genISH2}, and obtain  an unique  stochastic weak solution  $U_\eps  \in AC([0,T], D) \hat{\otimes}_\pi (S)_{-1} $ of the equation \eqref{SP reg}
whose chaos expansion is given by 
\begin{equation*}
U_\varepsilon(t, x, \omega) = \sum_{\gamma\in \mathcal I} u_{\varepsilon,\gamma}(t, x) \, H_\gamma(\omega), \quad t \in [0,T], \, x \in \R^d, \, \omega \in \Omega,
\end{equation*}
such that  
\begin{equation} \label{zvezda}
\normtri{U_\varepsilon}^2_{-p} < 3 M_\varepsilon(T)^2 A\left( 1 + 2C_{\frac{p}{2m}} C_{\frac{p(m-1)}{m}-s_\eps-2}\right), \quad p \geq  p_{U_\varepsilon},
\end{equation} 
where 
$$ p_{U_\varepsilon}:=  \big\lfloor \max\left\{2mp_F,2mp_G,\frac{m(3+s_\eps)}{m-1}\right\}\big\rfloor +1, \quad m \in \N,\, m>1, $$
 \begin{equation} \label{Seps}
 s_\eps := \left\{ \begin{array}{ll}
 	\frac{\ln{(\tilde{M}_\varepsilon(T) q_\varepsilon)^{2}}}{\ln 2} +1, & \text{for } (\tilde{M}_\varepsilon(T) q_\varepsilon)^2>1
 	\\
 	0, & \text{for } (\tilde{M}_\varepsilon(T) q_\varepsilon)^2 \leq 1,
 	\end{array}\right.
 \end{equation} 
and where $A$ defined in \eqref{A} is independent of $\eps$. 

In order to show that the net $(U_\varepsilon)_\varepsilon$ is $AC([0,T],L^2(\R^d)) \hat{\otimes}_\pi (S)_{-1}$-{moderate} in the sense of the Definition \ref{ModGSp2} we start from \eqref{zvezda}. 
The term $M_\varepsilon (T)$ is estimated as follows. Since it is possible to achieve the log-type moderateness of $q_{\mathbf{0}} \ast \varphi_{\varepsilon}$, see \cite[Remark 1.4]{BLO2022},  there exist $N_1 \in \N$ and $C_1>0$ independent of $\varepsilon$ such that 
\begin{equation} \label{ocenaM}
	M_\varepsilon(T) : = M \exp {\left( \left(w + M \|q_\mathbf{0}\ast \vphi_\eps \|_{L^\infty}\right)T \right)} \leq C_1 \eps^{-N_1},
\end{equation}
and thus
\begin{equation}\label{zvezda1}
\normtri{U_\varepsilon}^2_{-p} < 3 C_1^2\varepsilon^{-2N_1} A\left( 1 + 2C_{\frac{p}{2m}} 
C_{\frac{p(m-1)}{m}-s_\eps-2}\right),\quad p\geq p_{U_\eps}.
\end{equation}
From  \eqref{ocenaM} we have
\begin{equation} \label{ocenaTildaM}
	\tilde{M}_\varepsilon (T) = \int_0^T M_\varepsilon (s)\,ds \leq C_{11} \eps^{-N_1}, 
\end{equation}
where $C_{11} = C_1(T).$ From \eqref{star1} and from the  Assumption \ref{as_Qsing}, we have that 
$$ \|q_\gamma \ast \varphi_\eps\|_{L^\infty}  \leq 
  C_2 \eps^{N_2}, $$ 
with $N_2= s + d/2$ and  $C_2 = C_\varphi \lceil s\rceil  q$, with 
$q = \sup_{\gamma\in\mathcal{I}}  \|q_\gamma\|_{H^{-s}_c} <\infty$, which implies
\begin{equation} \label{qeps}
 q_\varepsilon : = \sup_{\gamma \in \mathcal{I}}\|q_{\gamma}\ast \vphi_\eps\|_{L^\infty} \leq C_2 \varepsilon^{-N_2}.
 \end{equation}
From \eqref{Seps}, using \eqref{ocenaTildaM} and \eqref{qeps}, one has 
\begin{equation} \label{ocenaSeps}
	s_{\eps} \leq \frac{2 \ln C \eps^{-N}}{\ln 2} + 1,
\end{equation} 
with $C=C_1C_{11}$ and $N=N_1+N_2,$ and $s_\eps \to \infty$ as $\eps \to 0.$
Choose 
\begin{equation}\label{peps}
p_{\eps}: = p_{U_\eps} + \frac{m}{m-1} \frac{1}{\eps} = \frac{m}{m-1}\left(s_\eps+3 + \frac{1}{\eps}\right), \quad \eps>0.
\end{equation} 
Then  
$p_{\eps} \geq  p_{U_\eps}$, for $\eps < \eps_0$, for some $\eps_0>0$,  thus \eqref{zvezda1} holds for $p\geq p_{\eps}$.  
 Using \eqref{ocenaSeps} one has that for all $\eps<\eps_0$ it holds
	$p_{\eps} \leq  2 \left(s_\eps +3+\frac{1}{\eps} \right) 
	\leq 2 \left(\frac{2\ln C\eps^{-N}}{\ln 2}  +  4 + \frac{1}{\eps} \right)$,
and hence, 
$$\lim_{\eps \to 0} p_{\eps} = \infty, \quad \text{ and } \quad
\lim_{\eps \to 0}  \left(\frac{p_{\varepsilon}(m-1)}{m}-s_\eps-2\right) =\lim_{\eps \to 0}  \left(1+\frac{1}{\eps}\right) = \infty.$$
Since  $C_p\to 0$ as $ p\to \infty$, 
we have that 
$$C_{\frac{p}{2m}} \to 0, \quad \text{ and } \quad C_{\frac{p(m-1)}{m}-s_\eps-2} \to 0,  \quad  p\geq p_{\eps},
\quad \text{as}\quad  \eps \to 0,$$
and therefore, for arbitrary  $\delta>0$, there exist $p_{0\eps} \in \N$ such that 
$C_{\frac{p}{2m}} C_{\frac{p(m-1)}{m}-s_\eps-2} < \delta,$ for all $p\geq p_{0\eps}$.
Thus, for $C=3C_1^2 A\left( 1 + \delta\right)$ and $N=2N_1$ independent of $\varepsilon$ $$\normtri{U_\varepsilon}_{-p} < C \varepsilon^{-N}, \quad p \geq \max\{p_{0\eps}, p_{\varepsilon} \},$$ 
and $(U_\eps)_\eps$ is moderate. Thus, $(U_\eps)_\eps$ is a stochastic very weak solution to \eqref{Eq:ISH3}.  \endproof

\section{Uniqueness of the stochastic very weak solution} \label{sec:singularQ_Uniqueness}
The question of the uniqueness of the stochastic very weak solution is understood via the notion of negligible  nets from the theory of generalized functions of Colombeau; see \cite{GKOS2001}. Two nets will be considered to represent the same solution if their difference is small enough, and tends to zero faster than any polynomial, as $\eps$ tends to zero. 
The uniqueness of the stochastic very weak solution holds  in the following sense.
\begin{thm}
	\label{T:unique} Let $s>0$, and let the potential {$Q\in H^{-s}_c(\mathbb R^d) \hat{\otimes}_\pi (S)_{-1}$} satisfy Assumption \ref{as_Qsing}.
	Let $(Q_\varepsilon)_{\varepsilon \in (0,1]}$ and $(\tilde{Q}_\varepsilon)_{\varepsilon \in (0,1]}$ be two different regularizing nets of  the singular potential $Q$ such that their difference is negligible, i.e., there exists $C>0$ such that for every $\varepsilon\in (0,1]$ and for all $n \in \mathbb{N}$
	\begin{equation} 
	\label{negligible_Q} 
	\normtri{ Q_\varepsilon -  \tilde{Q}_\varepsilon}_{-p} \leq  \, C \, \varepsilon^n, \quad  p \geq p_{Q}.
	\end{equation}
	Let  $(U_{\varepsilon})_{\varepsilon \in (0,1]} \subset AC([0,T];L^2(\R^d)) \hat{\otimes}_\pi (S)_{-1} $ and $( V_{\varepsilon})_{\varepsilon \in (0,1]} \subset AC([0,T];L^2(\R^d)) \hat{\otimes}_\pi (S)_{-1}$ be  two very weak solutions of  the stochastic  problem \eqref{Eq:ISH3} which correspond to  $(Q_\varepsilon)_{\varepsilon \in (0,1]}$ and  $(\tilde{Q}_\varepsilon)_{\varepsilon \in (0,1]}$, respectively. Then,  there exists $c>0$ such that
	for all $\varepsilon \in (0,1]$ and all $n\in \N$ it holds
	\begin{equation*}
	\normtri{U_\varepsilon - V_\varepsilon}_{-p} \leq c\, \varepsilon^n, \quad p\geq \max\{p_{U_\eps},p_{V_\eps},p_Q\}. 
	\end{equation*}
\end{thm}
\proof 
By assumptions, for each $\eps \in (0,1]$, the GSPs $U_\eps(t,x,\omega)  = \sum_{\gamma\in \mathcal I} u_{\varepsilon,\gamma}(t, x) \, H_\gamma(\omega)$, and $V_\eps (t,x,\omega)= \sum_{\gamma\in \mathcal I} v_{\varepsilon,\gamma}(t, x) \, H_\gamma(\omega)$ satisfy the regularized problem \eqref{SP reg} with potentials $Q_\eps=\sum_{\gamma \in \mathcal{I}} q_{\eps,\gamma}(x)H_\gamma (\omega)$, and $\tilde{Q}_\eps=\sum_{\gamma \in \mathcal{I}} \tilde{q}_{\eps,\gamma}(x)H_\gamma (\omega)$, respectively. Both nets $(U_\eps)_\eps$ and $(V_\eps)_\eps$ are $AC([0,T],L^2(\R^d))$-moderate, i.e., there exist $N_1,N_2 \in \N$ and $C_1,C_2 >0$ such that  
\begin{equation} \label{Ueps mod}
\normtri{U_\eps}_{-p} \leq C_1 \eps^{-N_1}, \, \normtri{V_\eps}_{-p} \leq C_2 \eps^{-N_2},\quad p\geq \max\{p_{\eps},\tilde{p}_{\eps}\}.
\end{equation}
%
%
 The GSP $U_\varepsilon-V_\varepsilon$ with the chaos expansion form
$$(U_\varepsilon-V_\varepsilon)(t, x, \omega) = \sum_{\gamma\in \mathcal I} (u_{\varepsilon,\gamma} - v_{\varepsilon,\gamma})(t, x) \, H_\gamma(\omega), \quad t \in [0,T], \, x \in \R^d, \, \omega \in \Omega,$$
satisfies the problem
$$\left(\dfrac{\partial}{\partial t} - \mathcal{L} \right) (U_\eps - V_\eps)+ Q_\eps \lozenge (U_\eps - V_\eps) = F_\eps, \quad (U_\eps - V_\eps) \big|_{t=0} = 0,$$
with $F_\eps = (\tilde{Q}_\eps - Q_\eps) \lozenge V_\eps \in AC([0,T],L^2(\R^d)) \hat{\otimes}_\pi (S)_{-1}.$ By \eqref{ocena} from Theorem \ref{T:genISH2}, we have that
$$\normtri{U_\eps-V_\eps}^2_{-p} \leq 3M_\eps (T)^2 T^2 \normtri{F_\eps}^2_{-p} \left(  1 + 2C_{\frac{p}{2m}} C_{\frac{p(m-1)}{m}-s_\eps-2}\right), \quad p \geq p_{U_\eps -V_\eps},$$
where $s_\eps$ is given by \eqref{Seps}.
Since $q_{\eps,\mathbf{0}} \ast \varphi_\eps$ is log-type moderate,  there exist $N_3 \in \N$ and $C_3>0$ independent of $\eps$ such that 
\begin{equation} \label{M mod}
M_{\eps}(T) \leq C_3 \eps^{-N_3}. 
\end{equation}
Since
$C_{\frac{p}{2m}}C_{\frac{p(m-1)}{m}-s_\eps-2} \to 0$ as $\eps \to 0$, there exists $\delta>0$ and $p_{0\eps} \in \N$ such that 
\begin{equation} \label{ocena Cp}
2C_{\frac{p}{2m}} C_{\frac{p(m-1)}{m}-s_\eps-2}<\delta, \quad p \geq \max\{p_{0\eps}, p_{U_\eps -V_\eps}\}.
\end{equation}
Therefore, from \eqref{M mod} and \eqref{ocena Cp} we obtain
\begin{equation}\label{estUminVeps}
\normtri{U_\eps-V_\eps}^2_{-p} < 3 C_3^2 \eps^{-2N_3} T^2 \normtri{F_\eps}^2_{-p} \left(  1 + \delta\right),\quad p\geq \max\{p_{0\eps}, p_{U_\eps -V_\eps}\},
\end{equation}
and we are left to estimate $\normtri{F_\eps}^2_{-p}$. It holds
\begin{eqnarray*}
	\normtri{F_\eps}^2_{-p} &=& \normtri{(\tilde{Q}_\eps - Q_\eps) \lozenge V_\eps}^2_{-p} = \sum_{\gamma \in \mathcal{I}} (2\N)^{-p \gamma}
\bigg	\| \sum_{\alpha + \beta =  \gamma} (\tilde{q}_{\eps,\alpha} -q_{\eps,\alpha} ) v_{\eps,\beta} 
\bigg	\|^2_{ AC([0,T],L^2(\R^d))}  \\
		& = & \sum_{\gamma \in \mathcal{I}}  
		\bigg\| \sum_{\alpha \in \mathcal{I}, \alpha \leq \gamma}  (\tilde{q}_{\eps,\alpha} -q_{\eps,\alpha} ) (2\N)^{-\frac{p \alpha}{2}}  v_{\eps,\gamma-\alpha}(2\N)^{-\frac{p (\gamma-\alpha)}{2}} 
		\bigg\|^2_{ AC([0,T],L^2(\R^d))} \\
		& = &  \sum_{\gamma \in \mathcal{I}}  \bigg(
		 \sup_{t\in[0,T]} 
		 \|  \sum_{\alpha \in \mathcal{I}, \alpha \leq \gamma}   (\tilde{q}_{\eps,\alpha} -q_{\eps,\alpha} ) (2\N)^{-\frac{p \alpha}{2}} v_{\eps,\gamma-\alpha} (t) (2\N)^{-\frac{p (\gamma-\alpha)}{2}} 
		 \|_{L^2(\R^d)}
		 \bigg)^2. 
\end{eqnarray*}
Using generalized H\"{o}lder inequality, we get
\begin{eqnarray*}	 
		\normtri{F_\eps}^2_{-p}	& \leq & 
 \sum_{\gamma \in \mathcal{I}} 
		\bigg( \sum_{\alpha  \in \mathcal{I}} \| \tilde{q}_{\eps,\alpha} -q_{\eps,\alpha} \|_{L^{\infty}(\R^d)}  (2\N)^{-\frac{p \alpha}{2}} 
		\bigg)^2
			\bigg( \sum_{\beta \in \mathcal{I}}  \sup_{t\in[0,T]}  \| v_{\eps,\beta}(t)\|_{L^2(\R^d)} (2\N)^{-\frac{p \beta}{2}}
		 \bigg)^2.
\end{eqnarray*}
Denoting by $B_1$ and $B_2$  numbers of terms in the  above sums with respect to $\alpha$ and $\beta$, respectively,  and applying \eqref{kvadratsume} 
we get
\begin{eqnarray*}
\normtri{F_\eps}^2_{-p}
		& \leq & \sum_{\gamma \in \mathcal{I}} 
		B_1 \bigg( \sum_{\alpha  \in \mathcal{I}} \| \tilde{q}_{\eps,\alpha} -q_{\eps,\alpha} \|_{L^{\infty}(\R^d)}^2  (2\N)^{-p \alpha} 
		\bigg) 	
		B_2 \bigg( \sum_{\beta \in \mathcal{I}}   \| v_{\eps,\beta}\|_{AC([0,T];L^2(\R^d))}^2 (2\N)^{-p \beta} \bigg)	\\
			&\leq&  B B_1 B_2 \normtri{\tilde{Q}_\eps - Q_{\eps}}^2_{-p}  \normtri{V_\eps}^2_{-p} , 
			\quad p\geq \max\{p_{Q}, p_{V_\eps}\},
\end{eqnarray*}		
where   $B$ is a constant denoting the number of terms in the above sum with respect to $\gamma$.

Therefore, using the above estimate on $\normtri{F_\eps}^2_{-p}$, negligibility of the difference of regularisations \eqref{negligible_Q}, and moderateness of the solutions  \eqref{Ueps mod} in \eqref{estUminVeps}
we obtain
\begin{equation*}
\normtri{U_\eps-V_\eps}^2_{-p} < 3 C_3^2 \eps^{-2N_3} T^2 B B_1 B_2 \normtri{\tilde{Q}_\eps - Q_\eps}^2_{-p} \normtri{V_\eps}^2_{-p} \left(  1 + \delta\right) =  c(T) \eps^{2n-2N_2-2N_3},
\end{equation*}
for arbitrary $n$ and $C = c(T) := 3 C_3^2 T^2 B B_1 B_2 C^2C_2^2 \left(  1 + \delta\right) $.
Choosing $n > N_2+N_3$ we obtain negligibility of the  difference of the solutions, i.e., 
$$\normtri{U_\eps-V_\eps}_{-p} \leq C \eps^k, \quad k \in \N, 
$$ 
with  $p\geq \max\{p_{\eps},\tilde{p}_{\eps}, p_{0\eps}, p_{U_\eps -V_\eps}, p_{Q}, p_{V_\eps}\}$, $\eps\to 0$. Therefore,  the stochastic very weak solution is unique in the quoted sence.
\endproof
%
%
%
%
%
\section{Consistency} \label{sec:singularQ_Consistency}
We aim to ensure that the different notions of solutions are consistent when all of them are well-defined. The consistency between the stochastic very weak solution to the stochastic parabolic problem \eqref{Eq:ISH3}, established in Theorem \ref{T:main}, and the stochastic weak solution from Theorem \ref{T:genISH2} is demonstrated in the following theorem.
\begin{thm} 
	\label{T:konz}  
Let the potential $Q\in L^\infty(\mathbb R^d) \hat{\otimes}_\pi (S)_{-1}$ satisfies Assumption \ref{as_Qbnd}, and let Assumptions \ref{op_L}, \ref{as_F} and \ref{as_G} hold. Let $V\in AC([0,T];L^2(\R^d)) \hat{\otimes}_\pi (S)_{-1}$ be the solution  to the problem \eqref{Eq:ISH3} obtained by Theorem \ref{T:genISH2}, and let $(U_\varepsilon)_{\varepsilon\in (0,1]} \subset AC([0,T];L^2(\R^d)) \hat{\otimes}_\pi (S)_{-1}$  be the very weak solution to the problem  \eqref{Eq:ISH3} obtained in Theorem \ref{T:main}.  
Then, for some $p>1$ 
	\begin{equation}
	\label{KonzKonver}
	\normtri{U_\varepsilon - V}_{-p} 
		 \to 0 \qquad \text{as} \qquad \varepsilon \to 0.
	\end{equation} 
\end{thm}
\proof
%
Since $(U_\varepsilon)_\varepsilon$ is a very weak solution to the problem \eqref{Eq:ISH3} obtained in Theorem \ref{T:main}, we have that for each $\varepsilon \in (0,1]$, $U_\varepsilon$ is a GSP in $AC([0,T];L^2(\R^d)) \hat{\otimes}_\pi (S)_{-1}$ with the chaos expansion 
$U_\varepsilon (t,x,\omega) = \sum_{\gamma \in \mathcal{I}} u_{\gamma,\varepsilon} (t,x) H_\gamma (\omega)$,
and it satisfies the regularized stochastic parabolic problem
\begin{equation} \label{eq:1}
	\left(\dfrac{\partial}{\partial t} -\mathcal{L} \right) U_\varepsilon (t,x,\omega) + Q_\varepsilon (x,\omega) \lozenge U_\eps (t,x,\omega) = F(t,x,\omega), \quad U_\eps(0,x,\omega) = G(x,\omega),
\end{equation}
where $Q_\varepsilon = Q \ast_x \varphi_\eps \in \D(\R^d) \hat{\otimes}_\pi (S)_{-1}$ is a regularisation of the potential $Q$.
Moreover, $(U_\varepsilon)_{\varepsilon \in (0,1]}$ is $AC([0,T];L^2(\R^d)) \hat{\otimes}_\pi (S)_{-1}$-moderate, i.e., there exist $N \in \mathbb{N}_0$ and $C>0$ such that for each $\eps \in (0,1]$ 
$$\normtri{U_\varepsilon}_{-p} \leq C \varepsilon^{-N},\quad p \geq p_\eps,$$ 
where $p_{\varepsilon}$ is defined in \eqref{peps}. The process $V \in AC([0,T];L^2(\R^d)) \hat{\otimes}_\pi (S)_{-1}$ is a stochastic weak solution to the problem \eqref{Eq:ISH3}, i.e., 
\begin{equation} \label{eq:2}
	\left(\dfrac{\partial}{\partial t} -\mathcal{L} \right) V(t,x,\omega) + Q (x,\omega) \lozenge V(t,x,\omega) = F(t,x,\omega), \quad V(0,x,\omega) = G(x,\omega),
\end{equation}
and it has the chaos expansion form 
$V(t,x,\omega) = \sum_{\gamma \in \mathcal{I}} v_\gamma (t,x) H_\gamma(\omega)$.
Subtracting \eqref{eq:1} and \eqref{eq:2},
adding and subtracting $Q_\varepsilon \lozenge V$, and then regrouping the factors, we obtain \begin{equation*} 
	\left(\dfrac{\partial}{\partial t} -\mathcal{L} \right) (U_\varepsilon - V) + Q_\varepsilon \lozenge (U_\varepsilon - V)  = (Q - Q_\eps) \lozenge V, \quad (U_\varepsilon - V)\big|_{t=0} = 0.
\end{equation*}
Therefore, for  all $\eps \in (0,1]$ the difference $U_\varepsilon - V$ satisfies the problem \eqref{Eq:ISH3}  in which a singular potential $Q$ is replaced by  a bounded potential $Q_\eps$, with zero initial condition, and load term $F_\eps := (Q - Q_\eps) \lozenge V \in AC([0,T];L^2(\R^d)) \hat{\otimes}_{\pi} (S)_{-1}$  with the chaos expansion form 
$$F_\eps (t,x,\omega) = \sum_{\gamma \in \mathcal{I}} \sum_{\alpha + \beta = \gamma} (q_\alpha - q_{\alpha,\eps})(x) v_{\beta} (t, x) H_\gamma (\omega).$$
Thus, by Theorem \ref{T:genISH2}, the difference $U_\varepsilon - V$ satisfies the estimate \eqref{ocena} which now takes the form
\begin{equation}\label{estUepsV}
	\normtri{U_\eps -V}^2_{-p} \leq 3M_\varepsilon (T)^2 T^2 \normtri{F_\eps}^2_{-p} \left(1+2C_{\frac{p}{2m}} C_{\frac{p(m-1)}{m}-s_\eps-2}\right), \quad p\geq p_{U_\eps -V},
\end{equation}
where $M_\eps (T)\! =\! M \exp((w+M\|q_{\mathbf{0},\eps}\|_{L^\infty})T)$, 
$ p_{U_\eps -V} \!:=\! \max \left\{ 2m  p_{F_\eps}, \frac{m(3+s_\eps)}{m-1} \right\}$, $m \in \N$, $m>1$,
 $s_\eps$ is given by \eqref{Seps}, and $q_\eps = \sup_{\gamma \in \mathcal{I}} \|q_{\gamma,\eps}\|_{L^\infty}<\infty.$ 
 By the Cauchy-Schwartz inequality, one has
 \begin{eqnarray*}
 	 \normtri{F_\eps}^2_{-p} & = & \sum_{\gamma 
 	 	\in \mathcal{I}} (2\N)^{-p\gamma}  \bigg\| \sum_{\alpha+\beta = \gamma} (q_\alpha - q_{\alpha,\eps}) v_\beta \bigg\|^2_{AC([0,T],L^2(\R^d))} \\
 	 	& \leq & \sum_{\gamma \in \mathcal{I}} B_{\alpha,\beta} \|(q_\alpha - q_{\alpha,\eps}) v_\beta\|^2_{AC([0,T],L^2(\R^d))} (2\N)^{-p\gamma}  \\
 	 	& = & \sum_{\gamma \in \mathcal{I}} B_{\alpha,\beta} \bigg(\|q_\alpha - q_{\alpha,\eps}\|_{L^\infty(\R^d) } \sup_{t \in [0,T]} \|v_\beta(t)\|_{L^2(\R^d)}\bigg)^2 (2\N)^{-p\gamma},
 \end{eqnarray*}
for $p\geq \max\{p_Q,p_V\}$, and  where $B_{\alpha,\beta}$ is the number of terms in the above sum.
 %
 As $\eps \to 0$,  we have that $\|q_\alpha - q_{\alpha,\eps}\|_{L^\infty(\R^d)} \to 0$,  $\normtri{F_\eps}^2_{-p} \to 0$, $M_\varepsilon (T) \to M(T)$, and 
 $C_{\frac{p}{2m}} C_{\frac{p(m-1)}{m}-s_\eps-2} \to 0$. 
%
 Therefore, from \eqref{estUepsV} we conclude  \eqref{KonzKonver}. \endproof

\section{Examples and concluding remarks} \label{sec6}
We consider an example of the problem \eqref{Eq:ISH3} taking the Laplasian for the operator $\mathcal{L}$, i.e.,
\begin{equation}
	\label{Example1}
		\begin{split}
\left(\frac{\partial}{\partial t}  - \Delta \right) \, U(t, x, \omega) + Q(x,\omega) \lozenge U(t, x, \omega) &= 
F(t,x,\omega), \enspace \\
			U(0, x, \omega) &= G(x,\omega),
		\end{split}
\end{equation}
where $t \in (0,T]$, $x \in \R$, $\omega \in \Omega.$ Note that the Laplasian is an unbounded closed operator with a dense domain $D=H_0^1(\R)$, and that it generates a (heat) $C_0$-semigroup of contractions $\{T_t\}_{t \geq 0}$ with stability constants $M=1$ and $w=0.$ For data we  assume the following.
\begin{itemize}
\item[$F:$] The force term $F\in AC([0,T];L^2(\R)) \hat{\otimes}_{\pi} (S)_{-1}$ is a Gaussian GSP\footnote{The stochastic process is a Gaussian process if and only if the nonzero
coefficients in its chaos expansion representation are only of the length 0 and
1, see \cite{LPS2015}}, 
\mbox{ defined by}
$$
F(t,x,\omega) = f(t,x) + g(x)W_{t}(\omega), \quad f \in AC([0,T],L^2(\R)), \enspace g\in L^2(\R),
$$
and $W_t$ is the time dependent white noise (see Example \ref{Ex:wn}). It has the chaos expansion
	$$F(t,x,\omega) =\sum_{\gamma \in \mathcal{I}} f_{\gamma}(t,x) H_\gamma (\omega) =f(t,x) + g(x) \sum_{k=1}^\infty \xi_k(t) H_{e_{k}} (\omega),$$
with $f_{\mathbf{0}}(t,x) = \mathbb{E}(F)=f(t,x)$,  $f_{e_k}(t,x)=g(x)\xi_k(t)$, $k\in\N$, and $f_\gamma(t,x)=0$, $\gamma \in \mathcal{I}\setminus \{\mathbf{0},e_k\}.$
Using that $\sup_{t \in [0,T]} |\xi_k (t)| =\mathcal{O}(k^{-\frac{1}{12}}),$ \cite[p. 22]{HOUZ1996}  for some constant $c>0$, one has
\begin{eqnarray*}
\normtri{F}^2_{-p} &=& \sum_{\gamma \in \mathcal{I}} \|f_\gamma (t,\cdot)\|^2_{AC([0,T];L^2(\R))} (2\N)^{-p\gamma} \\
& \leq &   \sup_{t \in [0.T]}  \|f(t,\cdot)\|^2_{L^2(\R)} + \|g(\cdot)\|^2_{L^2(\R)} \sum_{k=1}^\infty \|\xi_k (\cdot)\|^2_{L^\infty([0,T])} (2k)^{-p} \\
& \leq &  \sup_{t \in [0.T]}  \|f(t,\cdot)\|^2_{L^2(\R)} + \|g(\cdot)\|^2_{L^2(\R)} 2^{-p} c \sum_{k=1}^\infty k^{-\frac{1}{6}-p},
\end{eqnarray*}
which is finite for $p> \frac{5}{6}.$ Thus, the critical exponent for the process $F$ is $p_F =1.$

\item[$G$:]	 The initial condition $G=W_x \in \mathcal{S}(\R) \hat \otimes_{\pi} (S)_{-1} \subset L^2(\R) \hat \otimes_{\pi} (S)_{-1}$ is the space dependent white noise (see Example \ref{Ex:wn}) with zero expectation $\mathbb{E}(G)=0,$ and the chaos expansion 
	 $$W_x(\omega) =\sum_{\gamma \in \mathcal{I}} g_\gamma (x) H_\gamma(\omega) =\sum_{k=1}^\infty \xi_k(x) H_{e_k}(\omega),$$
where $g_{\mathbf{0}}(x)=0,$ $g_{e_k} (x) = \xi_k(x)$ for $k \in \N$, and $g_\gamma (x) =0$ for $\gamma \in \mathcal{I} \setminus \{\mathbf{0},e_k\}$. 
The sequence $\xi_k$, $k \in \N,$  forms an orthonormal basis of $L^2(\R)$, implying
$$\normtri{G}^2_{-p} = \sum_{\gamma \in \mathcal{I}} \|g_\gamma (\cdot)\|^2_{L^2(\R)} (2\N)^{-p\gamma} =\sum_{k=1}^{\infty} \|\xi_k(\cdot)\|^2_{L^2(\R)} (2k)^{-p} = \sum_{k=1}^\infty (2k)^{-p} < \infty,$$
for $p>1.$ Thus, the critical exponent for the process $G$ is $p_G =2.$

\item[$Q:$] For the singular potential $Q$ we take a Gaussian GSP, with the  chaos expansion 
	$$Q(x,\omega) =  \sum_{\gamma \in \mathcal{I}} q_{\gamma} (x) H_\gamma (\omega)=\delta(x-x_0) + \sum_{k=1}^{\infty} \delta (x - x_{e_k}) H_{e_k} (\omega), \quad x_0, \,x_{e_k} \in \R,$$ 	
where $\delta(x-a)$ denotes the Dirac delta distribution at point $x=a$, and coefficients are given as $q_{\mathbf{0}} (x) = \delta (x-x_0) = \mathbb{E}(Q),$ $q_{e_k} (x) = \delta (x-x_{e_k})$ for $k\in \N,$ and $q_{\gamma} (x) = 0$ for $\gamma \in \mathcal{I} \setminus \{\mathbf{0},e_k\}$. 
Thus, $Q\in  H_0^{-1}(\R)\hat \otimes_{\pi} (S)_{-1}$.  In addition, we have 
$$ q = \sup_{\gamma\in\mathcal{I}} \|\delta(\cdot-x_\gamma) \|_{H_0^{-1}} =  \sup_{\gamma\in\mathcal{I}} \sup_{\substack{ \|\phi\|_{H^1} =1}} | \langle\delta(x-x_\gamma),\phi  \rangle | = \sup_{\gamma\in\mathcal{I}} \sup_{\substack{\|\phi\|_{H^1} =1}} |\phi (x_\gamma)| = 1.$$
\end{itemize}
The singular random potential $Q$, 
the force term $F$, and the initial condition $G$  satisfy the assumptions of Theorem \ref{T:main}, and therefore, there exists a unique stochastic very weak solution $\net{U_\eps}$ to the  problem \eqref{Example1}. 
This means that for every $\eps >0$, $U_\eps$ is solution to the regularized problem 
\begin{equation}
	\label{RegEx}
	\begin{split}\small 
		\left( \frac{\partial}{\partial t}  - \Delta \right) \, U(t, x,\omega) + {Q_\varepsilon (x,\omega)} \lozenge  U(t, x,\omega) &= f(t, x) +g(x) W_t(\omega),	\\
		U(0, x,\omega) &= W_x(\omega),
	\end{split}
\end{equation}
with $\net{Q_\eps}$ being the  regularization of the singular  GSP $Q$ defined by \eqref{reg_Q}.  
 For later purposes we choose a  mollification net $\net{\varphi_\eps}$ defined as
$\varphi_{\varepsilon} (x):= \frac{1}{\lambda_\varepsilon} \, \varphi \big( \frac{x}{\lambda_\varepsilon}\big) \in \D(\mathbb{R})$, where $\lambda_\eps =\left( \log (\frac{1}{\eps}) \right)^{-1} $, and $\varphi$ is an arbitrary chosen mollifier.  
We obtain
$$Q_\eps(x,\omega) = \sum_{\gamma \in \mathcal{I}} q_{\eps,\gamma} (x) H_\gamma (\omega)=\varphi_\eps (x-x_0) + \sum_{k=1}^{\infty} \varphi_{\eps} (x - x_{e_k}) H_{e_k} (\omega), \quad \eps>0,$$
where  
the corresponding coefficients $q_{\eps, \gamma}$, $\gamma \in \{\mathbf{0},e_k\} $ satisfy
$\| \vphi_\eps (x-x_\gamma)\|_{L^\infty} =\frac{1}{\lambda_\eps} \|\varphi\|_{L^\infty} \leq C \log{\frac{1}{\eps}},$
 where 
 $$C = \|\varphi\|_{L^\infty} = \max_{x \in \supp \varphi} |\varphi (x)|< \infty.$$ This  insures that all coefficients in the chaos expansion of process $Q_\eps$ are log-type moderate.  In particular, for the zeroth order  coefficient, we have 
 \begin{equation}\label{ocena_q_eps0}
 \|q_{\eps,\mathbf{0}}\|_{L^\infty} = \|\varphi_\eps (\cdot-x_0)\|_{L^\infty}\leq C \log{\frac{1}{\eps}}.
 \end{equation}  
Furthermore, we have 
$$q_\eps: = \sup_{\gamma \in \{\mathbf{0},e_k\}}\| \vphi_\eps (\cdot-x_\gamma)\|_{L^\infty} =\frac{1}{\lambda_\eps} \|\varphi\|_{L^\infty},$$ 
which is finite for a fixed $\eps$.
The solution $U_\eps$ of the problem \eqref{RegEx} is given by
	\begin{equation} \label{ExSol}
U_\eps(t,x,\omega) = u_\mathbf{\eps,0}(t,x) + \sum_{k=1}^\infty u_{\eps,e_k}(t,x) H_{e_k} (\omega) + \sum_{|\gamma|>1} u_{\eps,\gamma} (t,x) H_\gamma (\omega),
\end{equation}
where the coefficients  $u_{\eps,\gamma}$, $\gamma\in\mathcal I$ are determined via
the chaos expansion method. We proceed to compute these coefficients explicitly following the
detailed exposition of the method given in the proof of Theorem 4 in \cite{GLO2023}.
Taking chaos expansions of all elements in \eqref{RegEx}, and collect corresponding terms to each $\gamma\in\mathcal{I}$, we have that
for $|\gamma|=0$ the expectation $\mathbb{E}(U_\eps) = u_\mathbf{\eps,0}$ satisfies the  problem
\begin{equation*}
	\begin{split}\small 
		\left( \frac{\partial}{\partial t}  - \Delta \right) \, u_\mathbf{\eps,0}(t,x) +\varphi_\eps (x-x_0) u_\mathbf{\eps,0}(t,x) &= f(t, x),	\\
		u_\mathbf{\eps,0}(0,x) &= 0,
	\end{split}
\end{equation*}
with a solution
 $$u_{\eps,\mathbf{0}} (t,x) = 
 \int_0^t S_{t-s} f(s,x)\,ds ,$$
 where $S_t = e^{t \varphi_\eps (x-x_0)} T_t$, $t \geq 0$,  is the semigroup generated by the perturbed Laplace operator $\Delta - \varphi_\eps (x-x_0) \text{Id}$, Id standing for the identity operator,
and for $|\gamma|=1$, we have that $\gamma =e_k$, $k \in \N$, and the coefficients $u_{\eps,e_k}$, $k \in \N,$ satisfy the deterministic PDEs
 \begin{equation*}
 	\begin{split}\small 
 		\left( \frac{\partial}{\partial t}  - \Delta \right) \, u_{\eps,e_k}(t,x) +\varphi_\eps (x-x_{0}) u_{\eps,e_k}(t,x) &= g(x)\xi_{k}(t)-\varphi_\eps (x-x_{e_k}) u_{\eps,\mathbf{0}}(t,x),	\\
 		u_{\eps,e_k}(0,x) &= \xi_k(x),
 	\end{split}
 \end{equation*}
with the solution
 $$u_{\eps,e_k}(t,x) = S_t \xi_k(x) + \int_0^t S_{t-s} (g(x)\xi_{k}(t)-\varphi_\eps (x-x_{e_k}) u_{\eps,\mathbf{0}}(t,x))\,ds.$$
 For  $|\gamma|>1$, the coefficients $u_{\eps,\gamma}$, $\eps>0,$ satisfy the regularized homogeneous problem \eqref{RegEx} with the initial condition $u_\mathbf{\eps,\gamma}(0,x) = 0$, implying that for $|\gamma|>1$, the coefficients $u_{\eps,\gamma} (t,x)$ vanish, i.e., $u_{\eps,\gamma} (t,x)=0$. As a result, the solution remains a Gaussian GSP, and the expression for the solution \eqref{ExSol} simplifies to
 $$U_\eps(t,x,\omega) = \int_0^t S_{t-s} f(s,x)\,ds + \sum_{k=1}^\infty \bigg( S_t \xi_k(x) + \int_0^t S_{t-s} (g(x)\xi_{k}(t)-\varphi_\eps (x-x_{e_k}) u_{\eps,\mathbf{0}}(t,x))\,ds\bigg) H_{e_k} (\omega).$$ 
 Moreover, using estimate \eqref{zvezda} from Theorem \ref{T:main}, and setting $m=2$, we obtain that for each $\eps>0$ it holds
 \begin{equation} \label{Ex:ocena}
 	\normtri{U_\varepsilon}^2_{-p} < 3 M_\varepsilon(T)^2 A\left( 1 + 2C_{\frac{p}{4}} C_{\frac{p}{2}-s_\eps-2}\right), \quad p \geq p_{U_\varepsilon},
 \end{equation}
where $ p_{U_\varepsilon}:=\max\left\{8, 2(3+s_\eps)\right\}$. By \eqref{ocenaM} and \eqref{ocena_q_eps0}, and stability constants for the semigroup $M=1$ and $w = 0$ we obtain that for $N<TC$ one has 
  $${M}_\eps(T)= e^{T\|q_{\eps,\mathbf{0}}\|_{L^\infty}} \leq e^{TC \log \frac{1}{\eps}} = \eps^{-TC} < \eps^{-N},$$  
  which, by the use of \eqref{mtilda}, leads to
  $$\tilde{M}_\varepsilon(T) =\frac{ e^{T\|q_{\eps,\mathbf{0}}\|_{L^\infty}}-1}{\|q_{\eps,\mathbf{0}}\|_{L^\infty} }
   = \frac{e^{\lambda_\eps^{-1} T\|\varphi\|_{L^\infty} }-1}{\lambda_\eps^{-1}\|\varphi\|_{L^\infty} }.$$ 
 According to \eqref{Seps} $s_\eps=0$  for  $\eps\geq \frac{T}{\ln 2}\|\varphi\|_{L^\infty}$, and for $\eps<\frac{T}{\ln 2}\|\varphi\|_{L^\infty}$  we have
    $$s_\eps= \frac{\ln{(\exp{(\lambda_\eps^{-1}T\|\varphi\|_{L^\infty}) }-1)^{2}}}{\ln 2} +1 \sim 
    \frac{2T\|\varphi\|_{L^\infty}}{\lambda_\eps \ln 2} +1.$$
 For arbitrary $\eta>0$ there exists $p_{0,\eps}$ such that $2C_{\frac{p}{4}} C_{\frac{p}{2}-s_\eps-2}\leq\eta$ for $p>p_{0,\eps}$.
Finally, from \eqref{Ex:ocena}, using the above calculations, for we derive the final estimate: there exists $N<TC$, with $C=\|\varphi\|_{L^\infty}$, and  $B=3  A\left( 1 + \eta\right)$ with the constant $A$ given by \eqref{A},
  $$	\normtri{U_\varepsilon}^2_{-p} < 3  A\left( 1 + \eta\right) \eps^{-2N}, \quad p \geq p_{U_\varepsilon},$$  
 yielding moderateness of the solution net $\net{U_\eps}$.

\subsection{Fields of applications}
The equation \eqref{Example1} appears in modeling various phenomenas, we mention several possible fields for applications.
\item[a)] {\em Disordered media with random impurities.} This field involves modeling the behavior of physical systems with structural irregularities, such as porous materials or composites.
The equation \eqref{Example1} models heat distribution in a medium, where $Q(x,\omega)$ represents the random distribution of highly localized impurities or defects in the medium which influence heat conduction or transport. In porous media, $Q(x,\omega)$  might represent spatially varying permeability influenced by random microstructures or fractures. One possible setting is described in \cite{FRANTZISKONIS}.

\item[b)] {\em Biological tissue with random heterogeneities.}  
Equation \eqref{Example1} models  biological processes, such as cell migration or diffusion, where the process $Q(x,\omega)$ is the term that accounts for influences of random variations in tissue density, random stiffness or compliance, or elastic properties, for example, while modeling drug delivery or tumor growth in heterogeneous tissues, or modelling articular cartilage with randomly distributed micro damage.  
%
As starting point one could take, for example, \cite{BP2000}.

\item[c)] {\em Anderson localization.} In quantum mechanics, equation \eqref{Example1} appears in the analysis of the behavior of electrons in random potentials. Homogeneous equation is used in the study of Anderson localization, where the random potential  $Q(x,\omega)$  describes spatially random energy barriers or wells affecting particle diffusion, see \cite{Anderson1958,LTW2009}.
\item[d)]  {\em Random risk potential in economy and finance.} 
The equation \eqref{Example1} is used in risk analysis to model random fluctuations in financial and economic systems. Applications include interest rate modeling, where 
$Q(x,\omega)$ represents spatially varying random risk factors, and option pricing under stochastic volatility, capturing the randomness in asset price behavior, see \cite{BS1973,ODTSB2024}.

\subsection{Further challenges} 
We outline several questions for future analysis:
\item[a)]
For future work,  it would be of significant interest to consider potentials $Q \in \E'(\R^d) \hat{\otimes}_{\pi} (S)_{-1}$,  with  chaos expansion $Q=\sum_{\gamma \in \mathcal{I}} q_\gamma H_\gamma,$ $q_\gamma \in H_c^{-s_\gamma}$,  $s_\gamma>0$, and  $\sum_{\gamma \in \mathcal{I}} \|q_\gamma\|_{H^{-s_\gamma}_c} (2\N)^{-p \gamma} < \infty$ for some $p \in \N$, see Remark \ref{Rem_on_reg}.

\item[b)]  Develop a framework for stochastic parabolic equations where the potential 
$Q$ is singular not only in space but also in time. 
Introduce an appropriate notion of regularization in both dimensions and analyze the resulting stochastic equations. 
%
Such models are relevant in systems where time-dependent impulses or shocks are present, such as neuronal dynamics or environmental transport with fluctuating sources.

\item[c)] Nonlinear extensions of \eqref{Eq:ISH3} also represent an intriguing avenue for research. 
For instance, these could include cases where the potential depends on the solution, such as 
$Q=Q(U)$, or where the force term becomes nonlinear, 
$F=F(U)$.
 For example,
$Q(U)= \Phi(U) +\eta$, where  $\Phi(U)$ introduces a feedback mechanism.
This introduces new challenges in defining and analyzing stochastic very weak solutions but opens up a broader range of applications in systems with feedback or self-organization. These investigations would broaden the scope of the current study and could pave the way for understanding more complex systems and applications.

\item[d)] Develop and analyze numerical schemes for solving \eqref{Eq:ISH3}  with singular potentials, employing adaptive mesh refinement for the spatial domain and stochastic Galerkin or Monte Carlo methods for the stochastic domain. Validate these solutions against analytical or experimental results, and conduct error analysis to quantify deviations between stochastic weak and very weak solutions. Establish convergence rates under Sobolev and Kondratiev norms, providing rigorous guarantees for the regularization approach and its applicability to complex systems.

 \item[f)]  Establish asymptotic properties of the solutions as $t\to\infty$. Show, for example, whether solutions decay to zero, stabilize, or exhibit oscillatory behavior depending on the nature of 
$Q=Q(U)$ and $F=F(U)$. This is critical in understanding the long-term behavior of stochastic heat or reaction-diffusion systems.

\item[g)] Investigate the effects of potentials $Q$ that exhibit long-range correlations in space or random
 variables. Derive conditions for the existence of solutions and study their qualitative behavior. This is particularly relevant in quantum systems with disordered potentials or biological systems with spatially distributed random heterogeneities.

\section*{Acknowledgments}

We are grateful to Professor G\"{u}nther H\"ormann from Faculty of Mathematics, University of Vienna, for discussions concerning constructions of appropriate Kondratiev spaces. 

\end{document}